\magnification=\magstep1
\input amstex
\documentstyle{amsppt}
\catcode`\@=11 \loadmathfont{rsfs}
\def\mycal{\mathfont@\rsfs}
\csname rsfs \endcsname \catcode`\@=\active

\vsize=6.5in

\topmatter 
\title Approximate equivalence of group actions \endtitle
\author  Andreas N\ae s Aaserud and Sorin Popa \endauthor

\rightheadtext{Approximate equivalence of actions}

\affil     {\it  University of California, Los Angeles} \endaffil

\address Math.Dept., UCLA, Los Angeles, CA 90095-1555\endaddress
\email  andreas.naes.aaserud\@gmail.com, popa\@math.ucla.edu\endemail

\thanks S.P. supported in part by NSF Grant DMS-1400208 \endthanks

\abstract  We consider several  weaker  versions of the notion of conjugacy  and orbit equivalence of measure preserving 
actions of countable groups on probability spaces, involving  equivalence of the ultrapower actions and asymptotic intertwining 
conditions. We 
compare them with the other existing  equivalence 
relations between group actions, and study the usual  type of 
rigidity questions around these new concepts (superrigidity, calculation of invariants, etc).  
\endabstract 

\endtopmatter

\document

\heading 1. Introduction \endheading

By a celebrated  result of Ornstein and Weiss ([OW1]), any two free ergodic probability measure preserving 
(pmp) actions of countable amenable groups $\Gamma \curvearrowright (X, \mu)$, $\Lambda \curvearrowright (Y, \nu)$ 
are orbit equivalent (OE). In turn, any non-amenable group $\Gamma$ is known to have ``many'' non-OE free ergodic pmp 
actions (cf. [CW], [H], [GP], [PS], [P6], [I1, GL, E]). Moreover, certain free ergodic pmp actions of non-amenable groups $\Gamma\curvearrowright X$ exhibit 
various degrees of rigidity, where mere orbit equivalence with another action $\Lambda \curvearrowright Y$ may imply isomorphism of 
the groups $\Gamma \simeq \Lambda$, and even isomorphism (conjugacy) of the two actions (see e.g. [Z],  [Fu1], [MS], [P7, P8, P9], [I2], [K1], [K2], etc). 
There has been constant interest in establishing such results, due to their intrinsic unexpected nature,  
their relevance to rigidity phenomena in von Neumann algebras and group theory,  
and their applications to Borel classification in descriptive set theory.   

There has also been much interest in obtaining existence of ``many'' distinct actions and rigidity results 
with respect to weaker notions of equivalence of actions, such as $W^*$-{\it equivalence} 
 (or von Neumann equivalence) considered in [P7], which  
requires the isomorphism of the associated group measure space von Neumann algebras $L^\infty(X)\rtimes \Gamma \simeq L^\infty(Y)\rtimes \Lambda$,  a condition well known to be weaker than OE of actions (cf. [S], [CJ], [PV]; 
for more on $W^*$-equivalence see [P10], [V2], [I3] and the references therein). More recently, 
a weak version of equivalence between actions of groups, 
that we will call here {\it weak conjugacy}, has been proposed in [Ke]. It    
requires that each one of the two actions can be ``simulated'' inside the other, in moments. 
Several results about weak conjugacy of actions have been obtained in [AE], [AW], [T], [IT], 
with many open problems still remaining.  

In this paper, we consider some additional weak versions of conjugacy and orbit equivalence of group actions, 
give examples and address the usual type of rigidity questions.  

All equivalences between actions that we study involve some form of approximation of the classical 
notions of conjugacy and orbit equivalence. As such, they can often be best formulated in the language of 
ultra products of actions and algebras, along some fixed (but arbitrary) free ultrafilter $\omega$ on $\Bbb N$.  
For instance,  using the ``ultraproduct framework'', weak conjugacy of pmp actions 
$\Gamma \curvearrowright^\sigma X$, $\Gamma \curvearrowright^\rho Y$ amounts to the 
ultrapower actions $\Gamma \curvearrowright^{\sigma^\omega} A^\omega$, $\Gamma \curvearrowright^{\rho^\omega}B^\omega$ 
(where $A=L^\infty(X)$, $B=L^\infty(Y)$)    
containing $\Gamma$-invariant subalgebras  $A_0\subset A^\omega$, $B_0\subset B^\omega$, 
such that $\Gamma \curvearrowright  A_0$ is isomorphic to $\Gamma \curvearrowright^\rho B$ and $\Gamma \curvearrowright B_0$ is isomorphic 
to   $\Gamma \curvearrowright^\sigma A$.

In this same spirit,  we will say that the two actions $\sigma, \rho$ are $\omega$-{\it conjugate} (respectively $\omega$ {\it orbit equivalent}, 
abbreviated $\omega$-{\it OE}), 
if the ultrapower actions $\Gamma \curvearrowright^{\sigma^\omega} A^\omega$, $\Lambda \curvearrowright^{\rho^\omega} B^\omega$ 
(respectively their full groups $[\sigma^\omega]$, $[\rho^\omega]$) are 
conjugate.  The actions $\sigma, \rho$ are {\it approximately conjugate} (respectively 
{\it approximately orbit equivalent}, abbreviated {\it app-OE}) if  they are $\omega$-conjugate (resp. $\omega$-OE) 
via an isomorphism $\theta: A^\omega \simeq B^\omega$ that's of ultrapower form, $\theta=(\theta_n)_n$, with $\theta_n : A \simeq B$, $\forall n$. 
This is easily seen to be equivalent to the existence of a sequence of isomorphisms $\theta_n: (X, \mu) \simeq (Y, \nu)$ 
that asymptotically intertwine the two actions (respectively the two orbit equivalence relations). 

The app-conjugacy of actions obviously implies $\omega$-conjugacy which implies weak conjugacy. 
We also consider an OE-version of weak conjugacy, which is  
weaker than $\omega$-OE (which in turn is implied 
by app-OE). We notice that the cost and the $L^2$-Betti numbers of an action $\Gamma \curvearrowright X$ (as defined in [Le], [G1] 
and respectively [G2]) 
are weak-OE invariant. We also show that co-rigidity ([A], [P3]) and relative Haagerup property ([Bo], [P5]) are weak-OE invariant.  
We prove that strong ergodicity of actions (as defined in [Sc]) is a  weak-OE invariant as well, and deduce from this and [CW] 
that  if a countable group $\Gamma$ is non-amenable and does not have property (T), then $\Gamma$ has at least two non 
weak-OE (thus non $\omega$-OE as well) free ergodic pmp actions (see Proposition 4.8 and Corollary 4.9).

We  show that if a group $\Gamma$ has an infinite amenable quotient (in particular, if $\Gamma$ has a free group as a quotient, 
or if it is the product between an infinite amenable group and an arbitrary group), then it has continuously many non-conjugate actions that are all approximately conjugate (see 3.7). 
On the other hand,  we notice that for a property (T) group $\Gamma$, app-conjugacy is the same as conjugacy  
(see 3.11). Related to this, it would be interesting  to decide whether any group that does not have property (T) admits  
many app-conjugate actions that are not conjugate. We also leave open the question of whether 
$\omega$-conjugacy is different from conjugacy for property (T) groups. 

By a result in [AW], all quotients of all Bernoulli actions are weakly equivalent, once they are free  (see also [T] and 5.1 in [P11] for more general results  
along these lines). 
Complementary to this, we show here that weak-conjugacy behaves well with respect 
to co-induction, so in particular if $\Gamma$ has an infinite amenable subgroup $H$, then the 
co-induction from $H$ to $\Gamma$ of any two free ergodic pmp actions of $H$ are weakly conjugate (Proposition 3.8).

Finally, we prove that if a strongly ergodic 
pmp action $\Gamma \curvearrowright X$ is OE-superrigid (i.e., any OE between $\Gamma \curvearrowright X$ 
and another free action $\Lambda \curvearrowright Y$ comes from a conjugacy), then it is app-OE superrigid (resp. $\omega$-OE superrigid), i.e., 
any app-OE (resp. $\omega$-OE) of this action with another free pmp action  $\Lambda \curvearrowright Y$ comes from an approximate conjugacy 
(resp. $\omega$-conjugacy) 
of $\Gamma \curvearrowright X$, $\Lambda \curvearrowright Y$ (see Proposition 5.1). We deduce from this (in fact, from the proof 
of this result) and the 
OE-superrigidity results in [K1, K2], [CK] that if 
$\Gamma$ is either a product of 
mapping class groups arising from compact orientable surfaces with higher complexity [K1],  or a certain type of amalgamated free product of higher rank lattices   
(such as $SL(3,\Bbb Z)*_\Sigma SL(3,\Bbb Z)$, with $\Sigma$ the subgroup of  matrices $(t_{ij})$ with $t_{31}=t_{32}=0$; see [K2]), 
or a central quotient of a surface braid group (cf. [CK]), then any   Bernoulli $\Gamma$-action (more generally, any strongly ergodic aperiodic $\Gamma$-action) 
is app-OE superrigid and $\omega$-OE superrigid 
(see Corollary 5.2).  We derive in a similar way an app-OE/$\omega$-OE strong rigidity result from OE strong rigidity  in [MS] (see Corollary 5.3). 

We have included several comments and open problems in the last part of Section 5. For instance,  given a free ergodic pmp action $\Gamma \curvearrowright X$ 
and a free ultrafilter $\omega$  on $\Bbb N$, 
we consider the 1-cohomology group of the Cartan inclusion $A^\omega \subset A^\omega \rtimes_{\sigma^\omega} \Gamma$, where $A=L^\infty(X)$ and 
$\sigma^\omega$ is the ultrapower action. This group, which we denote 
$H^1_\omega(\sigma)$, is an $\omega$-OE invariant for $\Gamma \curvearrowright^\sigma X$ (thus also app-OE invariant). It contains the group $(H^1(\sigma))^\omega$, where 
$H^1(\sigma)$ is the 1-cohomology group of $\sigma$, and we conjecture that for certain property (T) groups $\Gamma$ 
one should have equality, and therefore $H^1_\omega(\sigma)=H^1(\sigma)$ whenever $H^1(\sigma)$ is finite 
(e.g., when $\sigma$  is Bernoulli, cf. [PS]). 
Showing this would amount to proving vanishing results for 
$\Cal U(A)$-valued approximate cocycles  for $\Gamma\curvearrowright ^\sigma A$, i.e. for sequences of 
maps $w^n:\Gamma \rightarrow \Cal U(A)$ 
satisfying $\lim_\omega \|w^n_g \sigma_g(w^n_h)-w^n_{gh}\|_2=0$, $\forall g, h\in \Gamma$. But 
the deformation-rigidity arguments  
used  so far to prove such results for cocycles of Bernoulli actions break when passing to approximate cocycles, where a new idea seems to be needed. 

Such calculations of the app-OE invariant $H^1_\omega(\sigma)$ for Bernoulli 
actions and their quotients, as in ([P6]), would show in  particular that $\Gamma$ has   
infinitely many non $\omega$-OE actions that are all weak-OE (by [AW]).  We in fact leave 
open the problem of whether any property (T) group has at least two non $\omega$-OE free ergodic pmp actions. This question would of course be answered 
if one could show that any property (T) group $\Gamma$ has at least two  non weak-OE group actions. Recall in this respect that it is not even known  whether any 
property (T) group  $\Gamma$ 
has at least two non weakly conjugate actions (cf. [Ke], [T]).

{\it Acknowledgement}. We are  grateful to Lewis Bowen, Damien Gaboriau and Adrian Ioana  for many useful comments on a preliminary draft of this paper.

\heading 2. Basics on Cartan inclusions  \endheading 

The most basic example of operator algebras (and historically the first to be considered) were the group measure space 
von Neumann algebras, introduced by Murray and von Neumann in [MvN]. These are crossed-product von Neumann algebras 
arising  from groups acting freely on measure spaces, 
$\Gamma \curvearrowright (X,\mu)$, or equivalently from groups acting on function spaces, $\Gamma \curvearrowright L^\infty(X, \mu)$, and denoted $L^\infty(X) \rtimes \Gamma$. 
 Cartan subalgebras and Cartan inclusions of von Neumann algebras are abstractions of the inclusion $L^\infty(X)\subset L^\infty(X)\rtimes \Gamma$, 
as developed in [D], [S], [Dy], [FM]. We will recall in this section the definition,  
basic properties and main invariants of these objects.

Classically, these considerations are made in a separable framework, meaning that all von Neumann algebras 
are assumed to act on a separable Hilbert space, corresponding to the fact that they arise from countable groups acting 
on standard probability measure spaces,  $\Gamma \curvearrowright (X, \mu)$. This can be done by either 
emphasizing the countable equivalence relation $\Cal R_{\Gamma \curvearrowright X}$, given by the orbits of the action, as in [S], [FM],  
or by emphasizing the full group $[\Gamma]$ associated with the action, as in ([Dy], [P2], or 1.3 in [P6]). 
The notions of approximate equivalence of 
actions that we are interested in in this paper lead us to also consider ultrapowers of such group actions, 
$\Gamma \curvearrowright L^\infty(X)^\omega$, and the Cartan inclusion of non-separable von Neumann algebras  
$L^\infty(X)^\omega \subset L^\infty(X)^\omega \rtimes \Gamma$  that they entail. We  thus need to consider these objects  without any separability 
assumption, a fact that imposes the ``full group formalism''. However, all proofs are the same as in the separable case and will thus be 
only sketched, or even omitted.

\vskip .1in

\noindent
{\bf 2.1. Definition.} $1^\circ$ Let $(M, \tau)$ be a finite von Neumann algebra with a faithful normal trace state. A {\it Cartan subalgebra of $M$} 
is a maximal abelian $^*$subalgebra (MASA) $A\subset M$ whose {\it normalizer} in $M$, $\Cal N_M(A)=\{u\in \Cal U(M) \mid uAu^*A\}$, generates $M$ as a von Neumann algebra. 
We should note that in most cases considered, we will deal with Cartan subalgebras $A\subset M$ 
for which $\Cal N_M(A)$ is countably generated over $A$, in the sense that there exists a  
countable subgroup $\Cal N_0\subset \Cal N_M(A)$ such that $\Cal N_0\vee A=M$. Note that this is equivalent to $M$ having a countable 
orthonormal basis over $A$ (in the sense of [PP]), i.e. $\text{\rm dim}_A L^2(M)=\aleph_0$.  We also denote by $q\Cal N_M(A)$ the {\it quasi-normalizer} 
of $A$ in $M$, i.e. the set of partial isometries $v$ in $M$ with $vv^*, v^*v\in A$, $vAv^*=Avv^*$. It is easy to see that $q\Cal N_M(A)=\{uq \mid 
u\in \Cal N_M(A), q\in \Cal P(A)\}$ and that $q\Cal N_M(A)$ with the multiplication operation inherited from $M$ is a pseudogroup.  Note that one can always find a  
[PP]-orthonormal basis of $M$ over $A$ with elements in $q\Cal N_M(A)$ (see 2.1 in [P2]). 
Moreover, by (2.3-2.5 in [P2]), if $M$ is a factor, then 
there exists an orthonormal basis of $M$ over $A$ with unitary elements in $\Cal N_M(A)$.

\vskip .05in

$2^\circ$ Let  $(A, \tau)$ be an abelian von Neumann algebra with a normal faithful trace state. If $\Gamma$ is a subgroup of Aut$(A, \tau)$, 
then we denote by $[\Gamma]$ the {\it full group} generated by $\Gamma$, i.e. the group of all automorphisms 
$\theta\in \text{\rm Aut}(A,\tau)$ with the property that there exists a partition of $1$ with projections $p_n \in A$ and automorphisms $\theta_n \in \Gamma$, 
such that $\theta_n(p_n)$ are mutually disjoint and $\theta(a)=\Sigma_n \theta_n(ap_n)$, $\forall a\in A$.

We also denote by $[[\Gamma]]$ the {\it full pseudogroup} 
generated by $\Gamma$, i.e., the set of all partial isomorphisms $\phi: Ap \rightarrow Aq$, with $p, q\in \Cal P(A)$, for which there exists a 
partition of $p$ with projections $p_n\in A$ and automorphisms $\theta_n\in \Gamma$ such that $\Sigma_n \theta_n(p_n)=q$ and $\theta(a)=\Sigma_n \theta_n(ap_n)$, 
$\forall a\in Ap$. We denote by $[[\Gamma]]_0$ the sub-pseudogroup of identity partial isomorphisms in $[[\Gamma]]$, which coincides with the units of $[[\Gamma]]$, 
and which can be identified with $\Cal P(A)$. 
If $\Gamma \curvearrowright (A, \tau)$ is an action of $\Gamma$ by automorphisms then we still denote by $[[\Gamma]]$ (resp. $[\Gamma]$) 
the full pseudogroup (resp. full group) generated by its image in Aut$(A, \tau)$. 

Note that, 
by a standard maximality argument, we have $[[\Gamma]]=\{ \phi \mid \exists \theta \in [\Gamma], p\in \Cal P(A)$ such that $\theta_{|Ap}=\phi\}$ (see e.g. [Dy]).  

\vskip .05in

$3^\circ$ If $\phi$ is a partial isomorphism of $A$ (for instance, from some full pseudogroup $[[\Gamma]]$), then we denote by $L(\phi), R(\phi)\in \Cal P(A)$ its 
{\it left} and respectively {\it right supports} and by $p_\phi\in \Cal P(A)$ 
the maximal projection under $R(\phi)$ on which $\phi$ acts as the identity. 
Note that it is equal to the maximal projection $p$ under $L(\phi)$ 
for which $p \phi=p$ and that $p_{\phi^{-1}}=p_\phi$. If $\psi$ is another partial isomorphism on $A$, 
then we define a distance between $\phi$ and $\psi$ by letting $\|\phi -\psi\|_2 = ((\tau(R(\phi))+\tau(R(\psi))-2\tau(p_{\phi^{-1}\psi}))^{1/2}$.  
We will see in 2.5-2.6 below that this norm corresponds to the Hilbert-norm associated with the canonical trace of a full pseudogroup, 
thus justifying the notation $\|  \  \|_2$. It is an easy exercise to show that $\|\phi - \psi\|_2$ is equivalent to $\sup \{\|\phi(aR(\phi))-\psi(aR(\psi))\|_2 \mid a\in (A)_1\}$.

\vskip .1in

Note that if $A\subset M$ is a Cartan subalgebra in a finite von Neumann algebra, then $\{\text{\rm Ad}(u) \mid u\in \Cal N_M(A)\}=[\Cal N_M(A)]$ is a full group. 
Note further that this group is naturally isomorphic to $\Cal N_M(A)/\Cal U(A)$ and that the associated 
full pseudogroup $[[\Cal N_M(A)]]$  can be naturally identified with $q\Cal N_M(A)/\Cal I(A)$, where $\Cal I(A)$ denotes the set of partial isometries in $A$. 

We will see shortly that, conversely, if $[\Gamma]$ is a full group on $(A, \tau)$, then there exists a Cartan inclusion $(A\subset M, \tau)$ such that $[\Cal N_M(A)]=[\Gamma]$  (cf. 2.5 below). 

\vskip .1in

\noindent
{\bf 2.2. Examples.}  $1^\circ$ The typical example of a Cartan inclusion is provided by the Murray-von Neumann 
{\it group measure space construction}, as follows. Let $\Gamma \curvearrowright (X, \mu)$ be a probability 
measure preserving (pmp) free action of a countable group $\Gamma$.  Let $A=L^\infty (X)$ be endowed with the trace $\tau(a)= \int a \text{\rm d}\mu$ 
and consider the action $\Gamma \curvearrowright (A, \tau)$ given by $g(a)(t)=a(g^{-1}(t))$, $t\in X$, $g\in \Gamma$, $a\in A=L^\infty(X)$. 
Denote by $M_0$ the $^*$-algebra generated by a copy of  
$A$ and a copy of the group $\Gamma$, denoted $\{u_g\}_{g\in \Gamma}$ (the group of {\it canonical unitaries}), 
satisfying $u_g a u_g^*=g(a)$, $\forall g\in \Gamma$, $a\in A$. Thus, the elements in $M_0$ are formal finite sums $x=\Sigma_g a_g u_g$, with $a_g\in A$, 
with product given by $(a_gu_g)(a_hu_h)=a_gg(a_h)u_{gh}$ and $^*$-operation by 
$(a_gu_g)^*=g^{-1}(a^*)u_{g^{-1}}$. $M_0$ comes endowed with a trace state $\tau$ extending the trace on $A$, 
defined by $\tau(\Sigma_g a_g u_g)=\tau(a_e)$, where $e\in \Gamma$ is the neutral element. The completion of $M_0$ in the 
Hilbert norm given by $\|x\|_2=\tau(x^*x)^{1/2}$, $x\in M_0$, identifies naturally with the Hilbert space $\oplus_g L^2(X)u_g\simeq L^2(X)\otimes \ell^2\Gamma$, 
on which $M_0$ acts by left multiplication. The group measure space von Neumann algebra associated with $\Gamma\curvearrowright X$, 
denoted $L^\infty(X)\rtimes \Gamma$, is by definition the weak closure 
of $M_0$ in this representation. The fact that the action $\Gamma \curvearrowright X$ is free 
is equivalent to the fact that the subalgebra $A\subset M$ is a MASA in $M$. 
Since $M$ is generated by the unitaries in $A$ and $\{u_g\}_g$, which all normalize $A$, 
it follows that if  $\Gamma \curvearrowright X$ is free then $A\subset M$ is a Cartan inclusion.

\vskip .05in

$2^\circ$ One can generalize the group measure space construction from the case of free pmp actions to the case of 
pmp actions of countable groups $\Gamma \curvearrowright (X, \mu)$ that are not necessarily free, and more generally 
to trace preserving actions of arbitrary groups $\Gamma$ on abelian von Neumann algebras $(A, \tau)$. 
In the case $\Gamma$ is countable 
and $A=L^\infty(X)$, with $\Gamma \curvearrowright L^\infty(X)$ arising from a pmp action $\Gamma \curvearrowright X$, 
this construction can be described as follows (cf. [FM]). 

Let  $\Cal R=\Cal R_{\Gamma \curvearrowright X} \overset {\text{\rm def}}\to 
 = \{(t,gt)\mid t\in X, g\in \Gamma\}$ be the {\it orbit equivalence relation}  generated 
 by the action (viewed up to null sets in the first variable). Let $m$ be   
the unique measure on $\Cal R$ satisfying
$m(\{(t,gt)\mid t\in X_0\})=\mu(X_0)$, $\forall X_0\subset X$ and
$g\in \Gamma$. Each element $\xi \in L^2(\Cal R, m)$ can be viewed as a matrix $\xi=(\xi(t,t'))_{t\sim t'}$, with $(\int |x(t,t')|^2 \text{\rm d} m)^{1/2} < \infty$.  
For each automorphism $\theta\in \text{\rm Aut}(X, \mu)$ implemented by some $g\in \Gamma$, 
one denotes by $u_\theta$ the matrix with $u_\theta(t,t')$ equal to $1$ if $t'=\theta(t)$ and $0$ otherwise. 
Any finite sum $x=\Sigma_\theta a_\theta u_\theta$ acts on $\xi \in L^2(\Cal R, m)$ by matrix multiplication, 
$x\xi(t,t')=\Sigma_{s\sim t} x(t,s)\xi(s,t')$, $\forall
(t,t')\in \Cal R$. The set of such elements is a $^*$-subalgebra $M_0=\Cal B(L^2(\Cal R, m))$ whose weak closure is 
the von Neumann algebra associated with $\Cal R$. This algebra, denoted $L(\Cal R)$,  has $A\simeq L^\infty(X)\subset L(\Cal R)$ as a Cartan subalgebra (when viewed
as the set of matrices $x$ supported on the diagonal $\{(t,t)\mid
t\in X\}\subset \Cal R$). 

This construction can alternatively be described by using the full group $[\Gamma]$, or the full pseudogroup $[[\Gamma]]$ (see e.g. Section 1.3 in [P6]), 
this approach having the advantage of working for arbitrary $\Gamma$ and $A$ (not necessarily separable). We will describe a 
more general version of this construction (involving also a 2-cocycle on $[\Gamma]$) 
in Example 2.5 and Proposition 2.6 below. 

\vskip .05in
$3^\circ$ If $(A\subset M, \tau)$ is a Cartan inclusion (e.g., as in $2.2.1^\circ$ or $2.2.2^\circ$  above) 
and $\omega$ is a free ultrafilter on $\Bbb N$, then denote by $A^\omega\subset M^\omega$ 
the corresponding ultrapower inclusion (cf. [W]; see also Section 1.2 in [P11]). Denote by $M(\omega)$ the von Neumann subalgebra of $M^\omega$ 
generated by $\Cal N_M(A)$ and $A^\omega$, $M(\omega)=\Cal N_M(A)\vee A^\omega$. Note that if $M=A\rtimes \Gamma$, for some free action 
$\Gamma \curvearrowright (A, \tau)$, then $\Gamma$ also induces a free action on $A^\omega$, by $g((a_n)_n) = (g(a_n))_n$, where $(a_n)_n \in A^\omega$, 
and $M(\omega)$ naturally identifies with $A^\omega \rtimes \Gamma$. 

\vskip .1in
\noindent
{\bf 2.3. Definition.} $1^\circ$ A 2-{\it cocycle} for the full group $\Cal G$ on $(A, \tau)$ is a map $v: \Cal G\times \Cal G \rightarrow \Cal U(A)$ satisfying the  
the following properties 
$$
v_{\theta,\phi}v_{\theta \phi,\psi}=\theta(v_{\phi,\psi})v_{\theta,\phi \psi}, \forall \theta, \phi, \psi \in \Cal G; \tag a
$$
$$
p_{\theta_1\theta_2^{-1}}v_{\theta_1, \psi}=p_{\theta_1\theta_2^{-1} }v_{\theta_2,\psi}, 
\psi(p_{\theta_1\theta_2^{-1}})v_{\psi,\theta_1}=\psi(p_{\theta_1\theta_2^{-1}})v_{\psi,\theta_2}, 
\forall \theta_1, \theta_2, \psi\in \Cal G.  
$$

Such a cocycle is {\it normalized} if the following conditions hold true: 
$$
p_\theta v_{\theta,\psi}=p_\theta; \theta(p_\psi)v_{\theta,\psi}=\theta(p_\psi); 
v_{\theta, \psi}p_{\theta \psi}=p_{\theta \psi}, \forall \theta, \psi\in \Cal G. \tag b
$$

Two cocycles $v, v'$ are {\it equivalent}, $v\sim v'$, if there exists $w:\Cal G \rightarrow \Cal U(A)$ satisfying 
the condition 
$$
p_{\theta_1\theta_2^{-1}}w_{\theta_1}= p_{\theta_1\theta_2^{-1}}w_{\theta_2}, \forall  \theta_1, \theta_2 \in \Cal G, \tag c 
$$
such that 
$$
v'_{\theta,\psi} = w_\theta \theta(w_\psi) v_{\theta, \psi} w_{\theta\psi}^*, \forall \theta, \psi \in \Cal G. \tag d
$$  

It is easy to check that any $2$-cocycle is equivalent to a normalized $2$-cocycle. Moreover, one can check that if $\Gamma \subset \Cal G$ is a 
subgroup such that $[\Gamma]=\Cal G$ and one has a map $v:\Gamma \times \Gamma \rightarrow \Cal U(A)$ 
satisfying the properties $(a), (b)$ for all $\theta_1, \theta_2, \psi\in \Gamma$, then it extends uniquely to a normalized $2$-cocycle on $\Cal G$. 

Note that the set $Z^2(\Cal G)$ of cocycles on $\Cal G$ is 
an abelian group under multiplication as $\Cal U(A)$-valued functions on $\Cal G \times \Cal G$, and that the equivalence of 
2-cocycles amounts to equivalence modulo its subgroup $B^2(\Cal G)$ of 2-cocycles of the form 
$(\theta, \psi) \mapsto w_\theta \theta(w_\psi) w_{\theta\psi}^*$, with $w:\Cal G \rightarrow \Cal U(A)$ satisfying 
$p_{\theta_1\theta_2^{-1}}w_{\theta_1}= p_{\theta_1\theta_2^{-1}}w_{\theta_2}$, $\forall  \theta_1, \theta_2 \in \Cal G$.

\vskip .05in
$2^\circ$ A 2-cocycle for a full pseudogroup $[[\Cal G]]$ on $(A, \tau)$ is a map $v$ from $[[\Cal G]] \times [[\Cal G]]$ 
into $\Cal I(A)$ (the set of partial 
isometries in $A$) which formally satisfies exactly the same properties as $(a)$ above (but with $\theta, \theta_1, \theta_2, \phi, \psi$ in $[[\Cal G]]$). 
The  cocycle $v$ is normalized if the conditions $(b)$ are formally satisfied, $\forall \theta, \psi \in [[\Cal G]]$.  
Two normalized cocycles $v, v'$ for $[[\Cal G]]$ are equivalent if there exists 
$w:[[\Cal G]] \rightarrow \Cal I(A)$ such that  $(c)$ and $(d)$ are satisfied, $\forall \theta, \theta_1, \theta_2, \psi \in [[\Cal G]]$.

\vskip .1in
\noindent
{\bf 2.4. Example.} Let $\Gamma$ be a countable group and $\Gamma \curvearrowright^\sigma X$ a free ergodic pmp 
$\Gamma$-action. A map $v:\Gamma \times \Gamma \rightarrow \Cal U(A)$ satisfying the conditions  
$v_{g,h}v_{gh,k}=\sigma_g(v_{h,k})v_{g,hk}$, $\forall g,h,k\in \Gamma$, is called a $2$-cocycle of $\sigma$. 
We denote by $Z^2(\sigma)$ the multiplicative group of such 2-cocycles $v$ and call it the 
2-cohomology group of $\sigma$.  We also denote by $B^2(\sigma)$ its subgroup 
of 2-cocycles of the form $v_{g,h}=w_g \sigma_g(w_h) w_{gh}^*$, for some $w: \Gamma \rightarrow \Cal U(A)$.  

It is easy to check that any $v\in Z^2(\sigma)$ extends uniquely to  
a 2-cocycle $\tilde{v}$ of the full group $[\Gamma]$. If in addition $v$ satisfies 
$v_{e,h}=v_{g,e}=v_{g,g^{-1}}=1$, $\forall g,h\in \Gamma$,  then $\tilde{v}$ is a normalized 2-cocycle of $[\Gamma]$. 
The map $v\mapsto \tilde{v}$ implements an isomorphism  
$Z^2(\sigma)\simeq Z^2([\sigma])$ whose inverse is the restriction of $\tilde{v}$ to 
$\Gamma \times \Gamma$ and which takes $B^2(\sigma)$ onto $B^2([\sigma])$.  

Note that if we denote  
$Z^2(\Gamma, \Bbb T)=\{ \lambda: \Gamma \times \Gamma \rightarrow \Bbb T 
\mid \lambda_{g,h}\lambda_{gh,k}=\lambda_{h,k}\lambda_{g,hk},  \forall g,h,$ $k\in \Gamma \}$, 
the 2nd cohomology group of $\Gamma$, and by 
$B^2(\Gamma, \Bbb T)$ the subgroup of co-boundaries, then given any free pmp 
$\Gamma$-action $\Gamma \curvearrowright X$, any element 
$\lambda\in Z^2(\Gamma, \Bbb T)$ determines an element in $Z^2(\sigma)$, thus in $Z^2([\Gamma])$. We will discuss in 5.7.4 the 
problem of whether some $\lambda\in Z^2(\Gamma, \Bbb T)$ which is non-trivial in $H^2(\Gamma, \Bbb T)=Z^2(\Gamma, \Bbb T)/B^2(\Gamma,\Bbb T)$ is still 
 non-trivial in $Z^2([\Gamma])/B^2([\Gamma])$.

\vskip .1in
\noindent
{\bf 2.5. Example.} Let $A \subset M$ be a Cartan subalgebra.  We can then  write $[\Cal N_M(A)]$ as a well ordered set 
$\{\theta_i\}_{i\in I}$, with $id_A$ as its first element. We choose $u_{id}=1$. 
Assume that for the ``first'' $I_0\subset I$ elements in $I$, we have chosen elements $\{u_{\theta_i}\}_{i\in I_0}\subset \Cal N_M(A)$ 
such that if $i, j\in I_0$ and $q\in \Cal P(A)$ is so that $\theta_{i}, \theta_j$ agree on $Aq$, then $u_{\theta_i}q=u_{\theta_j}q$. 
If now $k$ is the first element of $I\setminus I_0$, then there is a maximal projection $p\in A$ such that the restriction of $\theta_k$ to $Ap$ 
does not agree with any $\theta_i, i\in I_0$, on $Aq$ for any $q\in \Cal P(Ap)$. Thus, we have 
$\theta_k= \oplus_{i\in I_0}{\theta_i}_{|Aq_i} \oplus {\theta_k}_{|Ap}$, for some mutually orthogonal projections $\{q_i\}_i$ and $p$. 
We define $u_{\theta_k}=\Sigma_{i\in I_0}u_{\theta_i}q_i + vp$, where $v\in \Cal N_M(A)$ is any unitary that implements $\theta_k$. 

In this manner, we have obtained a set of unitaries $\{u_\theta \mid \theta\in [\Cal N_M(A)]\} \subset \Cal N_M(A)$ such that $u_{id}=1$, 
Ad$(u_\theta)=\theta$, $\forall \theta$, and $u_\theta q =u_{\psi} q$ for any projection $q\in A$ with $\theta, \psi$ agreeing on $Aq$. 

Note that this implies  $p_{\theta_1\theta_2^{-1}}u_{\theta_1}u_{\theta_2}^{-1}=p_{\theta_1\theta_2^{-1}}$, $\forall \theta_1, \theta_2\in [\Cal N_M(A)]$. 
This in turn easily implies that $v_{\theta,\psi}=u_\theta u_\psi u^*_{\theta\psi}$, for $\theta, \psi\in  [\Cal N_M(A)]$,   
satisfies conditions $(a)$ in 2.3.1$^\circ$, and is 
thus a 2-cocycle for the full group $[\Cal N_M(A)]$. The choice of $v$ is unique up to 
the equivalence relation $\sim$ defined above. 

Thus, to each Cartan subalgebra $A\subset M$  
one can in fact associate a pair $(\Cal G, v/{\sim})$, consisting of a full group $\Cal G=[\Cal N_M(A)]$ on the abelian von Neumann algebra $(A, \tau)$ and 
the equivalence class of a $\Cal U(A)$ valued 
$2$-cocycle $v$ for $\Cal G$, as defined above. 
This provides a functor, from the category of tracial Cartan inclusions $(A\subset M, \tau)$ with morphisms given by 
trace preserving isomorphisms between the ambient algebras $(M, \tau)$ carrying 
the Cartan subalgebras $A$ onto each other, to the category of pairs $(\Cal G \curvearrowright (A, \tau), v/{\sim})$, 
consisting of a full group $\Cal G$ on an abelian von Neumann algebra $(A, \tau)$ and a (equivalence class of a) 2-cocycle 
$v:\Cal G \times \Cal G \rightarrow \Cal U(A)$, with morphisms given by trace preserving isomorphisms between the algebras $(A, \tau)$ 
that carry the full groups (resp. class of 2-cocycles) onto each other.  

This functor is in fact one to one and onto, its inverse being constructed as follows. Let $(\Cal G, v/{\sim})$ be a pair consisting 
of a full group $\Cal G$ on $(A, \tau)$ 
and a (class of a) normalized 2-cocycle $v:\Cal G \times \Cal G \rightarrow \Cal U(A)$ for 
$\Cal G\curvearrowright A$. Let $M_0$ be the vector space of finite sums $\Sigma_\theta a_\theta u_\theta$, with $a_\theta\in A$ 
and indeterminates $u_\theta, \theta\in \Cal G$ satisfying $u_\theta u_\phi = v_{\theta, \phi}u_{\theta\phi}$. On $M_0$ one defines 
the product rule $(a_\theta u_\theta) (a_\phi u_\phi)=a_\theta \theta(a_\phi)v_{\theta,\phi}u_{\theta\phi}$ 
and $^*$-operation given by $(a_\theta u_\theta)^*= \theta^{-1}(a_\theta^*)u_{\theta^{-1}}$. Define a functional $\tau$ on $M_0$ by  $\tau(a_\theta u_\theta)=\tau(a_\theta p_\theta)$, 
which clearly factors through the expectation $E$  of $M_0$ onto $A$ given by 
$E(\Sigma a_\theta u_\theta)=\Sigma a_\theta p_\theta$. Also, define 
the sesquilinear form on $M_0$ by $\langle x, y \rangle_\tau =\tau(y^*x)$, which is easily seen to be positive, semi-definite. 
Denote by $\Cal H$ the Hilbert  space completion of $M_0/\Cal I_\tau$, where $\Cal I_\tau = \{x\in M_0 \mid \tau(x^*x)=0\}$. Note that 
any $x\in M_0$ acts as a bounded linear operator on $\Cal H$,  by left multiplication. 

Finally, define $M=L(\Cal G,v)$  to be 
the weak closure of $M_0$ in this representation.  Then $\tau$ defines a normal faithful trace state on $L(\Cal G, v)$ and $A\subset L(\Cal G,v)$ is a Cartan inclusion. 
It is straightforward to check that if a pair $(\Cal G, v)$ comes  
from a Cartan inclusion $(A\subset M)$, as described above, then this new Cartan inclusion $A\subset L(\Cal G,v)$ naturally identifies 
with the initial one $(A\subset M)$; and that the pair $([\Cal N_M(A)], v)$ described above, for the Cartan inclusion $A\subset L(\Cal G, v)$, identifies naturally 
with $(\Cal G, v)$.

Altogether, we have just shown the following version of a well known result of Feldman-Moore ([FM]), 
formulated in terms of full groups rather than countable equivalence relations on a standard measure space, 
a fact that allows us avoid separability conditions.     

\proclaim{2.6. Proposition}  The functor from Cartan inclusions to pairs consisting of a full group and a normalized $2$-cocycle defined above, $(A\subset M, \tau) 
\mapsto ([\Cal N_M(A)], v/{\sim})$, is one to one and onto, its inverse being the functor 
$(\Cal G\curvearrowright (A,\tau), v/{\sim}) \mapsto (A\subset L(\Cal G,v))$. 
Via these functors, the isomorphisms of Cartan inclusions correspond to 
isomorphisms of full groups  intertwining the corresponding $($classes of$)$ $2$-cocycles. 
\endproclaim

We will often need to consider commuting squares of Cartan inclusions $A_0\subset M_0$, $A\subset M$, as also considered in (Section 1.1 of [P11]), 
under an additional non-degeneracy condition. Recall in this respect that two 
inclusions of finite von Neumann algebras  $B_0\subset N_0$, $B\subset N$,  
with $N_0\subset N$, $B_0\subset B$,  are in {\it commuting square} position if $B\cap N_0=B_0$ and the trace preserving expectations of 
$N$ onto $B$ and $N$ onto $N_0$ commute, with their product giving the $\tau$-preserving expectation of $N$ onto $B_0$, i.e. $E^N_BE^N_{N_0}=E^N_{N_0}E^N_B=E^N_{B_0}$ (see 1.2 in [P1]). For these conditions to hold true, it is in fact sufficient that $E^N_B(N_0)\subset N_0$ (equivalently $E^N_{N_0}(B)\subset B$). 

It is important to note that if $B_0\subset N_0$, $B\subset N$ are MASAs, with $B_0\subset B$, $N_0\subset N$, 
then the commuting square condition automatically holds true. 
This is because  $E^N_{N_0}(B)\subset B_0'\cap N_0=B_0$ (since $B_0\subset B$ are commutative).  

\vskip .1in
\noindent
{\bf 2.7. Definition.} Two Cartan subalgebras $A_0\subset M_0$,  
$A\subset M$  are said to form a {\it non-degenerate embedding} if $A_0\subset A, M_0\subset M$ (so by the above discussion, the commuting square relation $E_{A}(M_0)=A_0$ 
is automatically satisfied), $\Cal N_{M_0}(A_0)\subset \Cal N_M(A)$  and $\overline{\text{\rm sp}(M_0A)} = M$ (equivalently, $\Cal N_{M_0}(A_0)\vee A = M$). 
If this is the case, 
then we also say that $A_0\subset M_0$ is  {\it sub-Cartan inclusion} of $A\subset M$, or that $A\subset M$ is an {\it extension of} 
$A_0\subset M_0$. 

Note that for Cartan inclusions arising from full (pseudo)groups $\Cal G_0$ on $(A_0, \tau_0)$ and $\Cal G$ on $(A,\tau)$, with trivial 2-cocycle, 
such a Cartan embedding amounts to the existence of a subgroup $\Cal G_0'\subset \Cal G$ and a $\Cal G_0'$-invariant subalgebra $A_0'\subset A$, 
such that $\Cal G_0'$ generates $\Cal G$ as a full (pseudo)group on $A$  and $\Cal G_0\curvearrowright A_0$ is isomorphic to  $\Cal G_0'\curvearrowright A_0'$. 

In turn, if the Cartan inclusions involve separable von Neumann algebras, and we view them as coming from countable equivalence relations 
$\Cal R_0$ on $(X_0, \mu_0)$ and $\Cal R$ on $(X, \mu)$, then a non-degenerate Cartan  embedding of $L^\infty(X_0)\subset L(\Cal R_0)$ into $L^\infty(X)\subset L(\Cal R)$ 
amounts to a {\it local OE} (or {\it local isomorphism}) of $\Cal R, \Cal R_0$ in the sense of (Definition 1.4.2 in [P8]), i.e., an a.e. surjective measure preserving 
map $\Delta: (X,\mu) \rightarrow (X_0,\mu_0)$ for which there exists a measure $0$ subset $S\subset X$ such that for any $t\in X\setminus S$, 
$\Delta$ is a bijection from the $\Cal R$-orbit of $t$ onto the $\Cal R_0$-orbit of $\Delta(t)$. The terminologies (class bijective) {\it extension} 
and {\it orbit bijection} are also used for a map $\Delta$ satisfying such a property. 

Note that if $\Cal R=\Cal R_{\Gamma \curvearrowright X}$, $\Cal R_0=\Cal R_{\Gamma \curvearrowright X_0}$ 
for some pmp actions of a countable group $\Gamma$, and $\Delta: X \rightarrow X_0$ is a local isomorphism between $\Cal R, \Cal R_0$, then 
$\Gamma \curvearrowright X$ is free iff $\Gamma \curvearrowright X_0$ is free. If this is the case, then such a local-isomorphism corresponds to an 
actual quotient of the free action $\Gamma \curvearrowright X$ to a free action $\Gamma \curvearrowright X_0$. So indeed, one can view $\Cal R$ 
as an ``extension of $\Cal R_0$'', while $\Cal R_0$ can be viewed as a (free) ``quotient of $\Cal R$'' (when interpreting them as groupoids). Altogether, 
non-degenerate embeddings between Cartan inclusions can be seen as algebraic abstractions of  local-OE (or of quotient) between equivalence relations.  

\proclaim{2.8. Lemma} $1^\circ$ If a Cartan inclusion $A\subset M$ is such that $\text{\rm dim}_AM$ is countable, then for any $\Cal X\subset A$ countable subset and 
$\Cal U_0\subset \Cal N_M(A)$ 
countable subgroup with $\Cal U_0\vee A=M$, there exists an inclusion of separable von Neumann algebras $A_0\subset M_0$ 
with $\Cal X\subset A_0\subset A$, $M_0\subset M$, such that: $A_0$ is Cartan  in $M_0$;  $\Cal U_0\subset 
\Cal N_{M_0}(A_0)\subset \Cal N_M(A)$;  $A_0\subset M_0$ is a sub - Cartan inclusion of $A\subset M$. 

$2^\circ$ If $[\Gamma]$ is a full group on $(A, \tau)$ with $\Gamma$ countable and $X\subset A$ is a countable subset, then there exists a $\Gamma$-invariant 
separable von Neumann subalgebra $A_0\subset A$ containing $X$ 
and a subgroup $G_0\subset [\Gamma]$ that contains $\Gamma$ and leaves $A_0$ invariant,  
such that the full group of $\Gamma_{|A_0}$ on $\text{\rm Aut}(A_0, \tau)$ coincides with ${G_0}_{|A_0}$ and such that for any $\theta\in G_0$, the projection  
$p_\theta$ $($defined for $\theta$ viewed as an element in $[\Gamma])$, belongs to $A_0$. 
\endproclaim
\noindent
{\it Proof}. The two parts $1^\circ$ and $2^\circ$ are clearly equivalent by Proposition 2.6. The proof of $1^\circ$ can be easily adapted 
from the proof of (Section 1.2 Lemma in [P11], or Lemma 3.8 in [P12]), 
and is thus left as an exercise for the reader.  
\hfill $\square$

\vskip .1in
\noindent
{\bf 2.9. Definition.} A Cartan inclusion $A\subset M$ is {\it strongly ergodic} if for any $\delta_0 > 0$, 
there exist a finite set of unitaries $F\subset \Cal N_M(A)$ and $\varepsilon_0 > 0$, such that any $x\in A$, $\|x\|\leq 1$, that satisfies 
$\|[u, x]\|_2 \leq \varepsilon_0$, $\forall u\in F$, must satisfy $\|x-\tau(x)1\|_2\leq \delta_0$. 

If $M$ has countable orthonormal basis over $A$ (equivalently, $\Cal N_M(A)$ contains a countable subgroup $\Cal U_0$ such that $\Cal U_0 \vee A=M$), 
then strong ergodicity is easily seen to be equivalent to $M'\cap A^\omega=\Bbb C$, for some free ultrafilter $\omega$ on $\Bbb N$, and also equivalent to 
this being true for any free ultrafilter $\omega$. Moreover,  $A\subset M$ is not strongly ergodic 
if and only if $\forall F\subset \Cal N_M(A)$, $\varepsilon >0$, $\exists v\in \Cal U(A)$ 
such that $\tau(v)=0$ and $\|[u,v]\|_2 \leq \varepsilon$, $\forall u\in F$. 

A Cartan inclusion $A\subset M$ has 
{\it spectral gap} if for any $\delta_0 > 0$, 
there exist a finite set of unitaries $F\subset \Cal N_M(A)$ and $\varepsilon_0 > 0$, such that if $\xi \in L^2(A)\ominus \Bbb C$, $\|\xi\|_2 \leq 1$, satisfies 
$\|[u, \xi]\|_2 \leq \varepsilon_0$, $\forall u\in F$, then $\|\xi\|_2\leq \delta_0$. This condition obviously implies strong ergodicity. It has been  shown  
by K. Schmidt in [Sc1] (inspired by arguments in [C1]) that if $A\subset M$ are separable, then in fact the two conditions are equivalent 
(i.e., strong ergodicity implies spectral gap as well). 
His argument is easily seen to work for arbitrary (not necessarily separable) Cartan inclusions.  

A full group $\Cal G$ on $(A,\tau)$ is {\it strongly ergodic} (respectively has {\it spectral gap}) if its corresponding Cartan inclusion is strongly ergodic 
(resp. has spectral gap). 

\proclaim{2.10.  Proposition} $1^\circ$ Assume $A\subset M$ is a sub - Cartan 
inclusion of $B\subset N$. If  $B\subset N$ is strongly ergodic then $A\subset M$ is strongly ergodic. 

$2^\circ$ Let $A\subset M$ be a Cartan inclusion and $\omega$ a free ultrafilter. 
As in $2.2.3^\circ$, denote $A^\omega \subset M(\omega)=A^\omega \vee \Cal N_M(A)$. Then $A\subset M$ is strongly ergodic if 
and only if $A^\omega \subset M(\omega)$ is strongly ergodic. 
\endproclaim
\noindent
{\it Proof}. $1^\circ$ If $A\subset M$ is not strongly ergodic, then for any $F\subset \Cal N_M(A)$, $\varepsilon >0$, there exists $v\in \Cal U(A)$ with 
$\tau(v)=0$ and $\|[u,v]\|_2 \leq \varepsilon$, $\forall u\in F$. 

Let now $v_1, ..., v_n \subset \Cal N_N(B)$  and $\delta>0$. We will prove that there exists $v\in \Cal U(B)$ such that $\tau(v)=0$ and $\|[v,v_i]\|_2 \leq \delta$, $\forall i$. 

To do this, note first that for any $\alpha>0$, 
there exist $u_1, ..., u_m \in \Cal N_M(A)$, partitions of $1$ with projections $\{p_{ij}\}_{1\leq j \leq m} \subset B$,  $ 1 \leq i \leq n,$ and unitary elements $w_{i,j}\in B$,  
$1\leq i \leq n, 1\leq j \leq m$,  such that $\tau(1-\Sigma^m_{j=1}p_{ij})\leq \alpha^2$, $\forall i$, 
and $p_{ij}v_i= w_{i,j}p_{ij}u_j$, $\forall 1\leq i \leq n$, $1\leq j \leq m$. 

If we now take $v \in \Cal U(A)$ to have 
trace $0$ and to satisfy $\|[u_j, v]\|_2 \leq \varepsilon$, for some $\varepsilon >0$ and all $j\leq m$, then $v$ commutes with $p_{ij}, w_{i,j}$ and we get the estimates 
$$
\|v_i v - vv_i \|_2  \leq \|\Sigma_{j=1}^m( p_{ij}v_i v - v p_{ij} v_i)\|_2 + 2 \alpha 
$$
$$
=\|\Sigma_{j=1}^m (w_{i,j}p_{ij}u_j v - v w_{i,j}p_{ij}u_j)\|_2+2\alpha$$
$$ = (\Sigma_{j=1}^m \|p_{ij} w_{i,j} [v, u_j]\|_2^2)^{1/2} + 2\alpha \leq m^{1/2} \varepsilon + 2 \alpha. 
$$
Thus, if we take $\varepsilon < \delta/2m^{1/2}$, $\alpha \leq \delta/4$, then we get $\|[v_i, u]\|_2 \leq \delta$, $\forall i$. 

\vskip .05in 
$2^\circ$ By the first part, if $A^\omega \subset M(\omega)$ is strongly ergodic then $A\subset M$ is strongly ergodic. If in turn 
$A^\omega \subset M(\omega)$ is not strongly ergodic, then noticing that $\Cal N_M(A)\subset \Cal N_{M(\omega)}(A^\omega)$ 
it follows that for any $F\subset \Cal N_M(A)$  finite  and $\varepsilon > 0$, 
there exists $v\in \Cal U(A^\omega)$ of trace $0$ such that $\|[u,v]\|_2 \leq \varepsilon/2$, $\forall u\in F$. Thus, 
if $v=(v_n)_n$ with $v_n \in \Cal U(A)$, then $\lim_{n \rightarrow \omega} \|[u,v_n]\|_2 \leq \varepsilon/2$, $\forall u\in F$, implying that for some $n$ large enough we have 
$\|[u,v_n]\|_2 \leq \varepsilon $, $\forall u\in F$. Thus, $A\subset M$ is not strongly ergodic. 
\hfill $\square$

\vskip .1in

Two important OE invariants for a pmp group action $\Gamma \curvearrowright X$ are the cost ([L], [G1]) and Gaboriau's $L^2$-Betti numbers  ([G2]). 
They can also be viewed as invariants for the associated full group $[\Gamma]$ or Cartan inclusion  
$A=L^\infty(X)\subset L^\infty(X)\rtimes \Gamma =M$, a fact that allows us to extend the definitions to 
Cartan inclusions $A\subset M$ where $A$ is not necessarily separable.

\vskip .05in
\noindent
{\bf 2.11. Definitions.} $1^\circ$ If $[\Gamma]$ is a full group on $(A,\tau)$ and a subset $\{\phi_i\}_i \subset [[\Gamma]]$ 
generates $[[\Gamma]]$ as a full pseudogroup, then its cost $c(\{\phi_i\})$ is defined as $\Sigma_i \tau(L(\phi_i))\in [0, \infty]$. The {\it cost of} 
$[[\Gamma]]$ is equal to the infimum over all costs $c(\{\phi\}_i)$ of generating subsets $\{\phi_i\}_i\subset [[\Gamma]]$. 
If $A\subset M$ is a Cartan inclusion, then its {\it cost} is by definition the cost of $[[\Cal N_M(A)]]$.

\vskip .05in
$2^\circ$ If $A\subset M$ is a separable Cartan inclusion, then the $n$th $L^2$-Betti number of $A\subset M$, denoted $\beta_n^{(2)}(A\subset M)$, is by definition 
the $n$th $L^2$-Betti number (as defined in [G2]) of the countable equivalence relation $\Cal R_{A\subset M}$. Notice that by [G2], 
if $A_0\subset M_0$ is a sub - Cartan inclusion of $A\subset M$, then $\beta_n^{(2)}(A\subset M)=\beta_n^{(2)}(A_0\subset M_0)$. (N.B.: this amounts to 
showing that if $\Cal R_0$ is a free quotient of $\Cal R$, then $\beta^{(2)}_n(\Cal R)=\beta^{(2)}_n(\Cal R_0)$; this is in fact 
only proved in [G2] in case $\Cal R, \Cal R_0$ come from free actions of groups, but the proof for free quotients of groupoids is exactly the same; see also the end of Remark 9.24 in [PSV].)  

Let now $A\subset M$ be a Cartan inclusion with countable orthonormal basis but $A$ not necessarily separable. 
By Lemma 2.8, there exist separable sub-Cartan inclusions $A_0\subset M_0$ 
of $A\subset M$ and any two separable sub-Cartan inclusions are included into a larger separable sub-Cartan inclusion of $A\subset M$. 
The  $n$th  $L^2$-{\it Betti number} of $A\subset M$, denoted $\beta_n^{(2)} (A\subset M)$,  
is by definition the $n$th $L^2$-Betti number of any of its separable sub-Cartan inclusions (which by the above remark does not depend on the 
choice of the separable sub - Cartan inclusion).

\proclaim{2.12.  Proposition} $1^\circ$  If $A_0\subset M_0$ is a sub-Cartan inclusion of $A\subset M$, then $c(A_0\subset M_0) 
\geq c(A\subset M)$. 

$2^\circ$ Let $A\subset M$ be a Cartan inclusion and denote $A^\omega \subset M(\omega)=A^\omega \vee \Cal N_M(A)$. Then 
$c(A\subset M)=c(A^\omega \subset M(\omega))$. 

$3^\circ$ If $A_0\subset M_0$ is a sub-Cartan inclusion of $A\subset M$, 
then they have the same $L^2$-Betti numbers. 
Moreover, for a Cartan inclusion coming from  a free action, $A \subset A \rtimes \Gamma=M$, the $L^2$-Betti numbers of $A\subset M$ coincide with those of $\Gamma$. 

\endproclaim

\noindent
{\it Proof}. Part $1^\circ$ is trivial by the definitions while part 3$^\circ$ is implicit in ([G2]; see 
also the discussion at the end of 2.1.2$^\circ$ above). 

By $1^\circ$, in order to prove $2^\circ$ we only need to show that 
 $c=c(A\subset M) \leq c(A^\omega \subset M(\omega))=c_\omega$.  If $c_\omega=\infty$ then there is nothing to prove, so we may assume $c_\omega<\infty$.  Let $F\subset q\Cal N_M(A)$ be an arbitrary finite set and $\varepsilon >0$. Let now 
 $\{v_i\}_{i\geq 1} \subset q\Cal N_{M(\omega)}(A^\omega)$ be a set of partial isometries that generate 
 $[[\Cal N_{M(\omega)}(A^\omega)]]$ as a full pseudo-group and satisfies $c(\{v_i\}_i) < c_\omega + \varepsilon/2$.  
 
Represent each  $v_i$ as $(v_{i,m})_m$ 
 with $v_{i,m}\in q\Cal N_M(A)$. Since each $w\in F$ is a constant sequence when viewed in $M(\omega)\subset M^\omega$ and since it 
 can be approximated arbitrarily well  by an appropriate finite sum of reductions by projections in $A^\omega$ of 
 products of  $\{v_i\}_i$  and their adjoints, 
 it follows that there exists $m$ and $k$ such that  $c(\{v_{i,m}\}_{i \leq k})< c_\omega + \varepsilon/2$ 
 while  each  element in $F$  is $\varepsilon/2|F|$-contained in the full pseudogroup generated by $\{v_{i,m}\}_{i\leq k}$. 
 Thus, by Definition 2.1.3$^\circ$ we can find a set  
$E\subset q\Cal N_M(A)$ that contains  $\{v_{i,m}\}_{i\leq k}$ such that $c(E)\leq c_\omega + \varepsilon$ and such that the full pseudogroup it generates 
on $A$ contains $F$. Since $F$ was arbitrary, this shows that $c \leq c_\omega + \varepsilon$. But $\varepsilon >0$ was arbitrary as well, 
showing that $c\leq c_\omega$.  
\hfill $\square$

\vskip .1in 
Recall  now two more properties for a Cartan inclusion $A\subset M$: 

\vskip .05in
$(a)$ {\it co-rigidity} considered in  [P3] and [A] (see also Sec. 9 in [P4], or 5.6 in [P5]), 
an alternative terminology being  
``$M$  {\it has property} (T) {\it relative to} $A$''; 
\vskip .05in
$(b)$ {\it relative Haagerup property}, considered in [Bo] and in Sec. 2 of [P5]. 

\vskip .1in 

We define the same properties for an orbit equivalence 
relation $\Cal R_{\Gamma\curvearrowright X}$, and more generally for a full group $[\Gamma]$ on an  abelian 
von Neumann algebra $(A,\tau)$, by requiring that the associated Cartan inclusion $A\subset L([\Gamma])$ has the corresponding property. 

\proclaim{2.13.  Proposition} Assume $A_0\subset M_0$ is a sub-Cartan inclusion of  $A\subset M$.  
Then $A\subset M$ is co-rigid $($respectively has relative Haagerup property$)$ if and only if $A_0\subset M_0$ has this same property. 

\endproclaim

\noindent
{\it Proof}.  First note that one can find a $\|  \ \|_2$-dense 
countable subset  $\Cal X$ in $\Cal N_{M_0}(A_0)$ with the property that each $x\in \Cal X$ has finite spectrum.  This  can 
be seen by first taking an arbitrary countable $\|  \ \|_2$-dense subset  $\Cal X_0\subset \Cal N_{M_0}(A_0)$, 
then approximating each $x_0\in \Cal X_0$ by a sequence of unitaries  in $ \Cal N_{M_0}(A_0)$ that have finite spectrum (by using Rokhlin lemma) 
and then taking $\Cal X$ to be the set of all the resulting unitaries. 
Now note that for each $u\in \Cal X$, any orthonormal basis $\{m^u_j\}_j$ of $\{u\}'\cap A_0 \subset \{u\}'\cap A$ is an orthonormal basis for $A_0\subset A$ as well. 

Finally, notice that if $\Phi:M_0\rightarrow M_0$ is an $A_0$-bimodular $\tau$-preserving completely positive (cp) unital map and we take 
$u\in \Cal N_{M_0}(A_0)$,  with $\theta =\text{\rm Ad}u$ denoting the automorphism of $A_0$ that it implements, then $\Phi(u)a=
\Phi(ua)=\Phi(\theta(a)u)=\theta(a)\Phi(u)$, $\forall a\in A_0$. Thus, $\Phi(u)=ua_u$ for some $a_u\in A_0$. In particular, it follows that 
if $u\in \Cal X$ and $\{m^u_j\}_j\subset A$ as above, then $[\Phi(u), m^u_j]=0$, $\forall j$.  

Altogether, this shows that the proof of (Lemma 9.2 in [P4]) applies for the non-degenerate commuting square embedding of $A_0\subset M_0$ into $A\subset M$, 
implying that any $A_0$-bimodular $\tau$-preserving completely positive (cp) map $\Phi:M_0\rightarrow M_0$ 
extends (uniquely) to an $A$-bimodular $\tau$-preserving cp map $\tilde{\Phi}$ on $M$. 

Thus, the proof of (9.3 in [P4]) applies in exactly the same way to show that if $A\subset M$ is co-rigid then $A_0\subset M_0$ is co-rigid. On the other hand, 
the argument in the proof of 9.4 in [P4] (in fact a much simplified version of it) shows that if $A_0\subset M_0$ is co-rigid, then $A\subset M$ is co-rigid. 

For the relative Haagerup property, the proof that $A\subset M$ will have this property once $A_0\subset M_0$ has it, 
is quite trivial, by noticing that if $\Phi:M_0 \rightarrow M_0$ is $A_0$-bimodular, $\tau$-preserving, 
cp map that's {\it compact relative to} $A_0$, then its extension to an $A$-bimodular cp map $\tilde{\Phi}$ on $M$ constructed above is 
compact relative to $A$. The opposite implication is trivial.  
\hfill $\square$

\vskip .1in
\noindent
{\bf 2.14. Remark.}  Note that in view of the remarks we made after Definition 2.7, Proposition 2.10.1$^\circ$ shows in particular that 
if  $\Cal R$, $\Cal R_0$ are orbit equivalence relations arising from pmp actions of countable groups on standard probability spaces 
and $\Cal R$ is an (class bijective) extension of $\Cal R_0$, 
then $\Cal R$  strongly ergodic implies $\Cal R_0$ is strongly ergodic, while Proposition 2.13 shows that 
$\Cal R$ has Haagerup property (resp. is co-rigid) if and only if $\Cal R_0$ does.

\vskip .1in 

\noindent
{\bf 2.15. Definition.}  We end this section by recalling the definitions of the three {\it symmetry groups} of a Cartan inclusion 
$A\subset M$ (respectively orbit equivalence relation $\Cal R_\Gamma$, resp. full group $[\Gamma]$): 
\vskip .05in
$1^\circ$ The {\it automorphism group of} $A\subset M$, denoted Aut$(A\subset M)$, is the group of automorphisms of $(M, \tau)$ that leave $A$ 
invariant, and which we endow with the point-$\| \ \|_2$ convergence. Note that in case $M$ is separable, 
this is a Polish group. The {\it outer automorphism group} Out$(A\subset M)$ is the quotient of Aut$(A\subset M)$ by 
its (normal) subgroup Int$(A\subset M)=\{\text{\rm Ad}(u) \mid u \in \Cal N_M(A)\}$. Similarly, for equivalence relations $\Cal R_\Gamma$, 
we denote by Aut$(\Cal R_\Gamma)$ the group 
of automorphisms of $(X, \mu)$ that normalize $\Cal R_\Gamma$ (i.e. take each $\Gamma$-orbit onto a $\Gamma$-orbit) and 
by Out$(\Cal R_\Gamma)$ its quotient by the subgroup Int$(\Cal R_\Gamma)$ of automorphisms that leave each $\Gamma$-orbit fixed. 
Similarly for Aut$([\Gamma])$, Out$([\Gamma])$. 

\vskip .05in
$2^\circ$ We denote by Aut$_0(A\subset M)$ the subgroup of Aut$(A\subset M)$ consisting of all automorphisms of $M$ that 
act trivially on $A$, and by Out$_0(A\subset M)$ its quotient by the (normal) subgroup of inner 
such automorphisms, i.e., Int$_0(A\subset M)=\{\text{\rm Ad}(v) \mid v \in \Cal U(A)\}$. As noticed in ([S] and Sec. 1.2-1.5 in [P6]), 
one has a natural identification between Aut$_0(A\subset M)$ and the group $Z^1([\Gamma])$ of $\Cal U(A)$-valued $1$-cocycles   
for $[\Gamma] \curvearrowright A$, where $[\Gamma]=\Cal N_M(A)/\Cal U(A)$. Also, via this identification, Int$_0(A\subset M)$ corresponds 
to the group of coboundary $1$-cocycles $B^1([\Gamma])$, and thus Out$_0(A\subset M)$ corresponds to the first cohomology group 
$H^1([\Gamma])=Z^1([\Gamma])/B^1([\Gamma])$, which is by definition the {\it first cohomology group} $H^1(A\subset M)$ of the 
Cartan inclusion $A\subset M$. 

As in (1.1 of [P6]), in the case $\Gamma \curvearrowright^\sigma X$ is a free pmp action and 
$A=L^\infty(X)\subset L^\infty(X)\rtimes_\sigma \Gamma = M$, then we also denote $Z^1(\sigma)=Z^1(A\subset M)$, 
$B^1(\sigma)=B^1(A\subset M)$, $H^1(\sigma)=Z^1(\sigma)/B^1(\sigma)=H^1(A\subset M)$. 
Same in the case $\Gamma \curvearrowright^\sigma (A, \tau)$ is a free action 
on an abelian von Neumann algebra which is not necessarily separable. Note that in this case, an element $c\in Z^1(\sigma)$ 
corresponds to a map $c:\Gamma \rightarrow \Cal U(A)$ satisfying $c_g\sigma_g(c_h)=c_{gh}$, $\forall g,h\in \Gamma$ 
and the topology on $Z^1(\sigma)$ inherited from Aut$(A\subset M)$ corresponds to point convergence in $\|\ \|_2$: a net 
$c_i \in Z^1(\sigma)$ converges to $c\in Z^1(\sigma)$ if $\lim_i \|c_{i,g}-c_g\|_2=0$, $\forall g\in \Gamma$. 

Recall from ([Sc1], or Sec. 1.5 in [P6]) 
that if $A\subset M$ is strongly ergodic (which in the case $A\subset M=A\rtimes_\sigma \Gamma$ 
amounts to $\sigma$ being strongly ergodic) then  Int$_0(A\subset M)$ is closed in Aut$_0(A\subset M)$, equivalently 
$B^1(A\subset M)$ is closed in $Z^1(A\subset M)$. By [Sc1], if $M$ is separable, then the converse is in fact true as well. 

\vskip .05in 
$3^\circ$ If $A\subset M$ is a Cartan inclusion, with $M$ a II$_1$ factor, and $t> 0$, then its $t$-{\it amplification} 
$(A\subset M)^t$ is defined as the (isomorphism class of the) Cartan inclusion 
$(A\otimes D_n)p \subset p (M\otimes M_{n \times n}(\Bbb C))p$, for some $n \geq t$, 
where $p \in A \otimes D_n$ is a projection of (normalized) trace  $\tau(p)=t/n$, with $D_n \subset M_{n\times n}(\Bbb C)$ 
the diagonal subalgebra.  The {\it fundamental group} $\Cal F(A\subset M)$ 
of $A\subset M$ is the set of all $t>0$ with the property that $(A\subset M)^t \simeq (A\subset M)$. Since $((A\subset M)^t)^s=(A\subset M)^{ts}$, 
we see that $\Cal F(A\subset M)$ is a multiplicative subgroup of $\Bbb R_+^*=(0, \infty)$. 

The definition of $t$-amplification of a Cartan inclusion $A\subset M$ actually makes sense whenever $M$ is of type II$_1$, but not necessarily a factor, 
by taking $(A\subset M)^t$ to be the isomorphism class of $(A\otimes D_n)p \subset p (M\otimes M_{n \times n}(\Bbb C))p$ , for some $n \geq t$, 
with $p\in \Cal P(A\otimes D_n)$ a projection of central trace equal to $(t/n) 1_M$ (we leave it as an exercise to check 
that this doesn't depend on $n\geq t$, nor on the choice of $p$, and that in this more general case we still have 
$((A\subset M)^t)^s=(A\subset M)^{st}$). Then one defines the fundamental group $\Cal F(A\subset M)$ the same way. 

One should note that, with the above definition of amplification, Out$((A\subset M)^t)$ naturally identifies with Out$(A\subset M)$ 
and Out$_0((A\subset M)^t)=H^1((A\subset M)^t)$ with Out$_0(A\subset M)=H^1(A\subset M)$. In particular, when $A\subset M$ 
comes from an ergodic  pmp action of a countable group $\Gamma \curvearrowright X$,  these invariants are 
stable orbit equivalence invariants of $\Gamma \curvearrowright X$. 

At the same time, by Gaboriau's results in [G2], the $L^2$-Betti numbers of $A\subset M$ satisfy the scaling formula 
$\beta^{(2)}_n((A\subset M)^t)=\beta^{(2)}_n(A\subset M)/t$ and thus, if one of the $L^2$-Betti numbers is  non-zero and $\neq \infty$, 
then no amplification by $t\neq 1$ of $(A\subset M)$ can be isomorphic to $(A\subset M)$, i.e. $\Cal F(A\subset M)=\{ 1\}$.

\heading 3. Approximate conjugacy \endheading

Recall that two pmp group actions $\Gamma \curvearrowright^\sigma X$, $\Lambda \curvearrowright^\rho Y$ are {\it conjugate} (or {\it isomorphic}) if there exist  
isomorphisms $\Delta: (X, \mu) \simeq (Y, \nu)$, $\delta:\Gamma \rightarrow \Lambda$, such that 
$\Delta(\sigma_g(t)) = \rho_{\delta(g)}(\Delta(t))$, $\forall_{ae} t \in X$, and $\forall g\in \Gamma$. Note that this is equivalent to the existence 
of an integral preserving  isomorphism $\Delta: L^\infty(X) \simeq L^\infty(Y)$ 
such that for any $a\in L^\infty(X)$ we have $\Delta(\sigma_g(a))= \rho_{\delta(g)}(\Delta(a))$. Let us also recall the weak 
version of this notion proposed by Kechris in [Ke],  which merely requires that each one of the actions can be ``simulated'' inside the other, in moments: 

\vskip .1in 
\noindent
{\bf 3.1. Definition [Ke].} A pmp action $\Gamma \curvearrowright^\sigma X$ is {\it weakly embeddable} into  a pmp action $\Lambda \curvearrowright^\rho Y$ 
with respect to some isomorphism $\delta:\Gamma \simeq \Lambda$, if  
for any finite $F \subset \Gamma$, any $p_1, ..., p_n \in \Cal P(L^\infty(X))$, any $K\geq 1$ and  
any $\varepsilon >0$,  
there exist $q_1, ..., q_n \in \Cal P(L^\infty(Y))$ such that for any $1\leq k \leq K$,  any choice of elements $g_1, ..., g_k \in F$ and $1\leq n_1, ..., n_k \leq n$,  
we have 
$$
|\tau(\Pi_{i=1}^k\sigma_{g_i}(p_{n_i})) -\tau(\Pi_{i=1}^k\rho_{\delta(g_i)}(q_{n_i}))| 
< \varepsilon  \tag 3.1.1$$ 
The two actions $\Gamma \curvearrowright^\sigma X$, $\Lambda \curvearrowright^\rho Y$
are {\it weakly conjugate} if each one of them can be embedded into the other, with respect to the same group isomorphism. 

\vskip .05in 

Note that in fact the original definition of weak conjugacy in [Ke] only requires condition 3.1 to be satisfied for moments of degree $\leq 2$, i.e. for $K=1, 2$, 
but by (Proposition 10.1 in [Ke]), this is sufficient for it to hold true for all $K$.  We have chosen to formulate weak conjugacy in the form 3.1, 
as an  ``approximation in moments'' (simulation) condition,  
because it is closer to the spirit of our paper. This formulation has also the advantage of   translating into equivalent 
formulations in ultrapower framework, as follows (cf. [CKT]):  

\proclaim{3.2.  Proposition} Let $\Gamma \curvearrowright^\sigma X$, $\Lambda \curvearrowright^\rho Y$ be pmp actions of countable groups. Denote 
$A=L^\infty(X), B=L^\infty(Y)$. The following properties are equivalent: 

$(a)$ The actions $\sigma, \rho$ are weakly conjugate. 

$(b)$  There exists a free ultrafilter $\omega$ on $\Bbb N$, a $\Gamma$-invariant 
subalgebra $B_0\subset  A^\omega$, $\Lambda$ invariant subalgebra $A_0 \subset B^\omega$ $($of the corresponding ultrapower of $\sigma, \rho)$, 
such that $\sigma$ is isomorphic to $\Lambda \curvearrowright A_0$ and $\rho$ to $\Gamma \curvearrowright B_0$, with respect to the same 
identification $\Gamma \simeq \Lambda$.

$(c)$ Property $(b)$ holds true for any free ultrafilter $\omega$ on $\Bbb N$. 

\endproclaim
{\it Proof}.  Condition $(b)$ is clearly just a reformulation of  3.1 above, and it does not depend on $\omega$, so $(a), (b), (c)$ are all equivalent. 
\hfill $\square$

\vskip .1in 

Let us point out that results from [AW], [T], combined with 
results in [B1], [PS], [P6],  give plenty of examples of non-conjugate actions that are weakly conjugate.

\proclaim{3.3.  Proposition ([AW], [T])} Let $\Gamma$ be  a countable group. 

$1^\circ$ Any two free ergodic quotients of  Bernoulli $\Gamma$-actions 
are weakly conjugate. More generally, if $\{H^1_j\}_j, \{H^2_k\}_k \subset \Gamma$ are countable families 
of amenable subgroups and for each $i=1,2$, 
$\Gamma \curvearrowright^{\sigma_i} X_i$, 
is a free ergodic $\Gamma$-action which can be realized as the quotient of a generalized Bernoulli action $\Gamma \curvearrowright \Pi_j [0,1]^{\Gamma/H^i_j}$, 
then $\sigma_1, \sigma_2$ are weakly conjugate. 

$2^\circ$ If $\{H_j\}_j \subset \Gamma$ is a family 
of amenable subgroups, then given any free ergodic pmp action $\Gamma \curvearrowright^\sigma X$ and any quotient $\rho$ of a 
generalized Bernoulli action $\Gamma \curvearrowright  \Pi_j [0,1]^{\Gamma/H_j}$, the diagonal $\Gamma$-action $\sigma \times \rho$ is weakly conjugate to $\sigma$. 
\endproclaim
\noindent
{\it Proof}. The case of (quotients of) Bernoulli actions is a result in [AW], while the more general case of generalized Bernoulli actions 
and part $2^\circ$ are results from [T]. A slightly more general result, proved directly in ultrapower framework, can be found in (5.1 of [P11]\footnote{Note that  
the papers [AW], [T] are unfortunately not cited in [P11], as at the time of writing this paper the second named author was unaware of this prior work.}) 
\hfill $\square$

\vskip .1in

Note that if a group $\Gamma$ is sofic, then by results of Bowen in [B1], Bernoulli $\Gamma$-actions with different 
entropy base give non-conjugate actions, while by 3.3 above, they are all weakly conjugate. Also, for certain classes 
of groups $\Gamma$ (e.g. $\Gamma$ having property T, or $\Gamma$ a product of a non-amenable and an infinite group), 
then by [PS], [P6] there are many free quotients of generalized Bernoulli $\Gamma$-actions that are not conjugate, while by 3.3 they 
are all weakly conjugate. 

We'll now consider two additional notions of equivalence for group actions, 
situated between conjugacy and weak conjugacy. The first one requires straight conjugacy of the ultrapower actions.    

\vskip .05in

\noindent
{\bf 3.4. Definition.} Let $\omega$ be a free ultrafilter on $\Bbb N$. Let 
$\Gamma \curvearrowright^\sigma X$, $\Lambda \curvearrowright ^\rho Y$ be pmp actions of countable groups and denote $A=L^\infty(X)$, $B=L^\infty(Y)$. 
The two actions are {\it $\omega$-conjugate} if there exists a trace preserving isomorphism $\theta:A^\omega \simeq B^\omega$ 
that intertwines the ultrapower actions  $\sigma^\omega, \rho^\omega$, up to some group isomorphism $\delta:\Gamma \rightarrow \Lambda$. 

While such $\omega$-conjugacy is clearly an equivalence relation between group actions, it is not clear whether or not it depends on the 
chosen ultrafilter $\omega$. However, if we impose the conjugacy between $\sigma^\omega$ and $\rho^\omega$ 
to be of ultraproduct form $\theta=(\theta_n)_n$, with $\theta_n:A \simeq B$, $\forall n$, then the corresponding relation no longer depends on $\omega$ 
and can in fact be described ``locally'', without using the ultrapowers of the actions.  

From this point on, 
it may be useful to recall from Definition 2.1 that 
if $\theta_1, \theta_2\in \text{\rm Aut}(X,\mu)$, then  $\|\theta_1 - \theta_2\|_2=(2-2 \tau(p_{\theta_1\theta_2^{-1}}))^{1/2}$. More generally, 
recalling that if $\phi$ is  a local isomorphism  on $(X, \mu)$ then $p_\phi$ denotes the projection in $L^\infty(X)$ 
on which it acts as the identity, then for local isomorphisms $\phi_1, \phi_2$ 
we have denoted $\|\phi_1-\phi_2\|_2=(\tau(p_{\phi_1\phi_1^{-1}}) + \tau(p_{\phi_2\phi_2^{-1}}) - 2 \tau(p_{\phi_1\phi_2^{-1}}))^{1/2}$. 
If $\Cal G$ is a full pseudogroup on $(X, \mu)$ containing  $\phi_i$ (for instance $[[\{\phi_i\}_i ]]$), then this norm corresponds to the Hilbert norm 
$\|u_{\phi_1} - u_{\phi_2}\|_2$ given by the canonical trace on $L(\Cal G)$  that we have considered in 2.5, 2.6.

\proclaim{3.5.  Proposition} Let $\Gamma \curvearrowright^\sigma X$, $\Lambda \curvearrowright^\rho Y$ be pmp actions of countable groups. Denote 
$A=L^\infty(X), B=L^\infty(Y)$. The following 
properties are equivalent: 

$(a)$ There exists $\delta:\Gamma \simeq \Lambda$ such that  
for any finite set $F\subset \Gamma$ and any $\varepsilon >0$, 
there exists an isomorphism $\theta:(X, \mu)\simeq (Y,\nu)$ satisfying $\|\theta \sigma_g \theta^{-1}- \rho_{\delta(g)}\|_2\leq \varepsilon$, 
$\forall g\in F$. 

$(b)$ There exists $\delta:\Gamma \simeq \Lambda$ and a sequence of isomorphisms $\theta_n: (X,\mu)\simeq (Y,\nu)$ 
such that $\lim_n \| \theta_n \sigma_g \theta_n^{-1} - \rho_{\delta(g)}\|_2=0$, 
$\forall g\in \Gamma$.

$(c)$  Given any free ultrafilter $\omega$ on $\Bbb N$,  
there exist a group isomorphism $\delta:\Gamma \rightarrow \Lambda$  and a trace preserving isomorphism $\theta:A^\omega  \simeq B^\omega$, 
of ultra product form $\theta=(\theta_n)_n$, where $\theta_n : A \simeq B$, $\forall n$,  implementing an isomorphism of the ultrapower 
actions $\Gamma \curvearrowright^{\sigma^\omega} A^\omega$, $\Gamma \curvearrowright^{\rho^\omega} B^\omega$, i.e., 
$\theta \sigma^\omega_g \theta^{-1}=\rho^\omega_{\delta(g)}$, $\forall g\in \Gamma$.

$(d)$ Property $(c)$ holds true for some free ultrafilter $\omega$ on $\Bbb N$. 
\endproclaim
\noindent
{\it Proof}. Clearly $(a) \Leftrightarrow (b) \Rightarrow (c) \Rightarrow (d)$. If $(d)$ holds and $\theta=(\theta_n)_n : A^\omega \simeq B^\omega$ is the isomorphism 
satisfying $\theta \sigma_g^\omega \theta^{-1} = \rho_{\delta(g)}^\omega$, $\forall g$, then 

$$
\lim_{n \rightarrow \omega}  \sup \{\| \theta_n(\sigma_g(a))  - \rho_{\delta(g)}(\theta_n(a))\|_2 \mid a\in (A)_1\} = 0, \forall g\in \Gamma \tag 3.5.1
$$ 

Let now $F \subset \Gamma$ be a finite set and $\varepsilon >0$. By  (3.5.1), there exists a neighborhood $\Cal V$ of $\omega \in \overline{\Bbb N}$ 
such that, when we view $\Cal V$ as a subset of $\Bbb N$, for any $m\in \Cal V$ we have: 

$$
\sup\{ \| \theta_{m}(\sigma_{g}(a))  - \rho_{\delta(g)}(\theta_{m}(a))\|_2  \mid a\in (A)_1 \} \leq \varepsilon, \forall  g \in F \tag 3.5.2
$$ 
showing that the actions $\sigma, \rho$ satisfy $(a)$.  
\hfill $\square$

\vskip .1in

\noindent
{\bf 3.6. Definition.} Two pmp actions of countable groups $\Gamma \curvearrowright^\sigma X$, $\Lambda \curvearrowright^\rho Y$ 
are  {\it approximately conjugate} ({\it app-conjugate}) if  any of the equivalent conditions in 3.5 holds true. 

\vskip .1in

Note that we obviously have the implications ``conjugacy  $\Rightarrow$ approximate conjugacy $\Rightarrow$ $\omega$-conjugacy $\Rightarrow$ weak conjugacy''. 
It is not clear whether there exist $\omega$-conjugate actions that are not app-conjugate. It would also be interesting to decide whether  the property of being 
$\omega$-conjugate for two group actions is independent of $\omega$.

If $\Gamma$ is a countable amenable group, 
then all free ergodic pmp $\Gamma$-actions are app-equivalent. This can be easily deduced from the Ornstein-Weiss' hyperfinite approximation of actions  
of amenable groups (or directly from the Rohlin lemma in [OW1]). It has in fact already been pointed out in (5.2 of [P11]). Indeed,  
if $\Gamma \curvearrowright^\sigma X$, $\Gamma \curvearrowright^\rho Y$ are two such actions, 
and $A=L^\infty(X)$, $B=L^\infty(Y)$, 
then  by  (5.2 in [P11]) there exists $v\in \Cal N_{M^\omega}(A^\omega)$ 
that intertwines the actions $\Gamma \curvearrowright^{\sigma^\omega} A^\omega$, $\Gamma \curvearrowright^{\rho^\omega} B^\omega=A^\omega$. 
But $v$ can be represented as a sequence $(v_n)_n$ with $v_n \in \Cal N_M(A)$, whence the automorphism Ad$(v)$ of $A^\omega$ is of ultraproduct type $(\text{\rm Ad}v_n)_n$. 

From the app-equivalence of $\Gamma$-actions in the case $\Gamma$ is amenable,  we will now derive that if a group $\Gamma$ has an infinite amenable quotient, 
then it has many non-conjugate actions that are approximately conjugate. Note that the class of such groups includes 
the free groups $\Gamma=\Bbb F_n, 2\leq n \leq \infty$, product groups $\Gamma = H \times K$, 
with $H$ infinite amenable and $K$ arbitrary, and groups of the form $\Gamma=\Bbb Z^2 \rtimes \Gamma_0$, with $\Gamma_0\subset SL(2,\Bbb Z)$.

\proclaim{3.7. Proposition} Let $\Gamma$ be a countable group which has an infinite amenable quotient $\pi:\Gamma \rightarrow H$. 
If $H \curvearrowright^{\rho_i} Y_i$, $i\in I$, are free ergodic pmp actions and $\Gamma \curvearrowright^\rho Y$ is a weak mixing pmp action 
which is free on $\ker \pi$, 
then the actions $\sigma_i$ of $\Gamma$ on $X_i=Y_i \times Y$ defined by $\sigma_i(g)(t,s)=(\rho_i(\pi(g))(t), \rho(g)(s))$ 
are free, ergodic and mutually app-conjugate $($thus $\omega$-conjugate, $\forall \omega$,  and weakly-conjugate as well$)$. 
If in addition $\rho_i, i\in I,$ are Bernoulli $H$-actions with distinct  entropy and $\rho$ is mixing on $\ker \pi$, 
then $\sigma_i, i\in I,$ are mutually non-conjugate. 
\endproclaim
\noindent
{\it Proof}. Since obviously the (diagonal) product of actions behaves well to app-conjugacy, the first part follows 
from the arguments above. 

If one takes two actions $\sigma_i, \sigma_j$ as above and $\delta: \Gamma \simeq \Gamma$ is an automorphism of $\Gamma$ 
such that $\sigma_i$ is conjugate with $\sigma_j$ with respect to $\delta$, 
and $\Delta: X_i \simeq X_j$ implements this conjugacy, then this is the same as $\sigma_i$,  $\sigma_j \circ \delta$ 
being intertwined by $\Delta$.  Now note that if we denote $K=\ker\pi$, then for any automorphism $\delta'$ of $\Gamma$, 
the quotient of $\Gamma$ by $\delta'(K)\cap K$ is still amenable. 
In particular, $K_0=\delta^{-1}(K)\cap K$ is infinite, with $K_0\subset K$, $\delta(K_0)\subset K$. Since $\Delta \sigma_i(K_0) \Delta^{-1}= \sigma_j(\delta(K_0))$, 
the fixed point algebra $(L^\infty(X_i))^{\sigma_i(K_0)}=:B_i$ is taken by $\Delta$ onto the fixed point algebra $(L^\infty(X_j))^{\sigma_j(\delta(K_0))}=:B_j$. 

Since $K_0$ is contained in $K$ and is infinite, it follows that $\sigma_i(K_0)$ is mixing on $L^\infty(Y)=1\otimes L^\infty(Y)\subset L^\infty(X_i)$ 
and acts trivially on $L^\infty(Y_i)=L^\infty(Y_i)\otimes 1$. Thus, $B_i=L^\infty(Y_i)$. In exactly the same way we deduce that $B_j=L^\infty(Y_j)$. 

Thus, $\Delta(L^\infty(Y_i))=L^\infty(Y_j)$. Since $\Delta$ intertwines $\sigma_i, \sigma_j \circ \delta$, this also implies that $\delta(\ker \pi)=\ker \pi$ 
and that $\Delta$ actually implements a conjugacy of the Bernoulli $H$-actions $\rho_i, \rho_j$. Since $\sigma_i, i\in I$, 
have distinct Kolmogorof-Sinai type entropies for different $i$'s (as defined for actions of 
amenable groups by Ornstein and Weiss in [OW2]), 
this  shows that if $\sigma_i, \sigma_j$ are conjugate, then $i=j$.   
\hfill $\square$ 

\vskip .1in

As we mentioned in Proposition 3.3, results in [AW] show that all free quotients 
of Bernoulli actions of a group $\Gamma$ are weakly conjugate. The previous result derives  
weak conjugacy of certain actions from their app-conjugacy. We note below that the co-induction of two weak-conjugate actions gives rise to   
weak conjugacy as well.  

To this end, let us first recall the construction of the co-induction of group actions (cf. [Lu]). Thus, 
let  $H\subset \Gamma$ be an inclusion of groups  
and $H\curvearrowright^\sigma (A,\tau)$ a trace preserving action on  an abelian 
von Neumann algebra with a faithful normal trace $\tau$ (but with $A$ not necessarily separable).  
Let $S\subset \Gamma$ be a set of representatives for $\Gamma/H$ and define $\pi: \Gamma \rightarrow H$ by $\pi(sh)=h$, $\forall s\in S, h\in H$. 
The map $c:\Gamma \times (\Gamma/H) \rightarrow H$ defined by $c(g, sH)=\pi(gs)\pi(s)^{-1}$, $s\in S, g\in \Gamma$, 
is then a 1-cocycle for the action $\Gamma \curvearrowright \Gamma/H$, 
with different choices of $\pi$ giving equivalent cocycles.  

Let $(\tilde{A}, \tau)=(A, \tau)^{\overline{\otimes} \Gamma/H}$. We define on it a trace preserving  $\Gamma$-action $\tilde{\sigma}$ as follows. 
For each $a\in A$ and $s\in S$ denote $\tilde{a}_{sH} = \cdots 1 \otimes a \otimes 1 \cdots \in A^{\otimes \Gamma/H}$, where $a$ appears on the $sH$ position. 
If $g\in \Gamma$ we then let $\tilde{\sigma}_g(\tilde{a}_{sH})= \tilde{b}_{gsH}$, where $b=\sigma_{c(g,sH)}(a)$. The cocycle relation implies that 
 $\tilde{\sigma}_{g_1} (\tilde{\sigma}_{g_2}(\tilde{a}_{sH}))=\tilde{\sigma}_{g_1g_2} (\tilde{a}_{sH})$, $\forall g_1, g_2\in \Gamma$, 
$s\in S$, $a\in A$, and thus it extends to an action $\tilde{\sigma}$ of $\Gamma$ on $\tilde{A}$. This is what is called the co-induction from $H$ to $\Gamma$ 
of $H \curvearrowright^\sigma A$. It is well known that if $\sigma$ is free then $\tilde{\sigma}$ is free.

Note that if $A_0\subset A$ is an $H$-invariant subalgebra, then one can 
naturally identify the action $\Gamma \curvearrowright\tilde{A_0}$ obtained as the co-induction from $H$ to $\Gamma$ of $H \curvearrowright A_0$ 
with the restriction of  $\tilde{\sigma}$ to the $\Gamma$-invariant subalgebra $\tilde{A_0}=A_0^{\otimes \Gamma/H}\subset A^{\otimes \Gamma/H}=\tilde{A}$. 

\proclaim{3.8.  Proposition} Let $\Gamma$ be a countable group with an infinite subgroup $H\subset \Gamma$. Let   
$H \curvearrowright^{\sigma_i} X_i$ be a free ergodic pmp action and denote $\Gamma \curvearrowright^{\tilde{\sigma_i}} \tilde{X_i}$ its   
corresponding  co-induction to $\Gamma$, $i=1,2$. If $H \curvearrowright^{\sigma_1} X_1$ is weakly-conjugate to 
$H \curvearrowright^{\sigma_2} X_2$, then $\Gamma \curvearrowright^{\tilde{\sigma_1}} \tilde{X_1}$ is weakly-conjugate   
to $\Gamma \curvearrowright^{\tilde{\sigma_2}} \tilde{X_2}$.   
\endproclaim
\noindent
{\it Proof}. This is trivial from the above observations, by taking into account that if $H\curvearrowright^{\rho_0} A_0$ is an action, 
then the co-induction  from $H$ to $\Gamma$ of its ultrapower $\rho_0^\omega$ identifies with the restriction to a $\Gamma$-invariant subalgebra of the ultrapower of 
$\Gamma \curvearrowright^{\tilde{\rho_0}}\tilde{A_0}$. 
\hfill $\square$

\vskip .1in 

Note that 3.8 gives an alternative proof of the weak-equivalence of Bernoulli $\Gamma$ actions ([AW]) for all groups $\Gamma$ 
that contain an infinite amenable subgroup (so a particular case of Proposition 3.3). More generally we have:

\proclaim{3.9.  Corollary} If a group $\Gamma$ has an infinite amenable subgroup $H\subset \Gamma$, then 
all pmp $\Gamma$-actions  that are co-induced from 
free ergodic pmp $H$-actions, are weakly-equivalent. 
\endproclaim
\noindent
{\it Proof}. We have already noticed that any two free ergodic pmp $H$-actions are app-conjugate. By Proposition 3.8 it follows that the co-induction 
from $H$ to $\Gamma$ of these $H$-actions are weakly-conjugate as well. 
\hfill $\square$

\vskip .1in

While Proposition 3.7 provides a rather large class 
of groups $\Gamma$ that have many non-conjugate free ergodic pmp actions which are all  app-conjugate, these groups  $\Gamma$ 
need to have infinite amenable quotients. Thus, they cannot have property (T). In fact, we notice below that for property (T)  groups app-conjugacy is the same 
as conjugacy. The proof reproduces arguments in ([P7], [P8], proof of 6.1 in [PV]; see also 14.2 in [Ke]), but we have included a full proof for 
completeness.

\proclaim{3.10.  Lemma}  Let $\Gamma$ be a property $(\text{\rm T})$ group and $F\subset \Gamma$ a finite set of generators. 
For any $\varepsilon >0$ there exists $\delta >0$ such that if $\Gamma \curvearrowright^\sigma (X, \mu)$, $\Gamma \curvearrowright^\rho (Y,\nu)$ are free 
ergodic pmp $\Gamma$-actions and $\theta: (X, \mu) \simeq (Y,\nu)$ satisfies $\|\theta \sigma_g \theta^{-1} - \rho_g\|_2 \leq \delta$, $\forall g \in F$, 
then there exists $\theta' : (X,\mu)\simeq (Y,\nu)$ such that $\theta' \sigma_g = \rho_g \theta'$, $\forall g\in \Gamma$, and $\|\theta'-\theta \|_2 \leq \varepsilon$. 
\endproclaim
\noindent
{\it Proof}. By property (T), we can choose $\delta>0$ such that, if $\pi$ is a unitary representation of $\Gamma$ on a Hilbert space $\Cal K$ 
and $\xi\in \Cal K$, $\|\xi\|=1$, satisfies $\|\pi_g(\xi)-\xi\|\leq \delta$, $\forall g\in F$, then the projection $\xi_0$ of $\xi$ 
onto the space of vectors fixed by $\Gamma$ satisfies $\|\xi-\xi_0\|\leq \varepsilon^2/4$. We will prove that this  $\delta$ checks the required condition. 

Let $A=L^\infty(X)$ with $\tau$ the trace on $A$ implemented by $\mu$. Denote by $\Cal G$ the full group 
generated in Aut$(A,\tau)$ by $\theta^{-1} \rho_\Gamma \theta $ and $\sigma_\Gamma$. 
Denote $A\subset L(\Cal G)=M$ the Cartan inclusion  
of this full group and for each $\phi\in \Cal G$, denote by $u_\phi$ the corresponding canonical unitary. 

Let $\langle M, e_A \rangle$  
denote as usual the basic construction von Neumann algebra for  $A\subset M$, generated in $\Cal B(L^2(M))$ 
by $M$ (viewed as left multiplication operators $L_x$ by elements $x\in M$) and by the orthogonal projection $e_A$, of $L^2(M)$ onto 
$L^2(A)$. Thus, $\langle M, e_A \rangle$ is the weak closure of sums of elements of the form $xe_Ay$, 
with $x, y\in M$, and it has a canonical faithful normal semifinite trace $Tr$ defined by $Tr(xe_Ay)=\tau(xy)$. Also, $xe_Ay \mapsto xy$ 
gives an $M$-bimodular operator valued normal semifinite faithful weight $\Phi$ of $\langle M, e_A \rangle$ 
onto $M$ satisfying $\tau(\Phi(xe_Ay))=Tr(xe_Ay)$.

The algebra $\langle M, e_A \rangle$ can also be described as the commutant of $JAJ$ in  $\Cal B(L^2M)$, where $J=J_M$ is the canonical 
involution on $L^2(M)$. It follows that if $u, v$ are unitaries in $M$ that normalize $A$, then 
Ad$(u) \text{\rm Ad}(JvJ) $ normalizes $\langle M,e_A \rangle$, acting on elements 
of the form $xe_Ay$ by Ad$(u) \text{\rm Ad}(JvJ) (xe_A y)=uxv^* e_A v yu^*$ and thus preserving the trace $Tr$. 

Moreover, since both $\text{\rm Ad}(u)$ and $\text{\rm Ad}(JvJ)$ 
leave invariant $A$ and $JAJ$, their product leaves invariant $\tilde{A}=A \vee JAJ = \vee_u u(Ae_A)u^*$. Finally, note that 
for any $u, v$ as above,  Ad$(u) \text{\rm Ad}(JvJ) $ leaves $M$ invariant, acting on it as Ad$(u)$, and that $\Phi(\text{\rm Ad}(u) \text{\rm Ad}(JvJ) (X))
=u\Phi(X)u^*,$ $\forall X \in \text{\rm sp} Me_AM$, $\Phi(\tilde{A})=A$. 

Denote $u_g=u_{\sigma_g}$ and $v_g=u_{\theta^{-1}\rho_g\theta}$, $g\in \Gamma$. Thus, we can define a unitary 
representation $\pi$ of $\Gamma$ on the Hilbert space  
$\Cal H=L^2(\langle M,e_A \rangle, Tr)$, by $\pi_g = \text{\rm Ad}(u_{g}) \text{\rm Ad}(Jv_gJ)$. 
Since $\pi_g(e_A)=u_gv_g^*e_Av_gu_g^*$ and $\|\sigma_g - \theta^{-1} \rho_g \theta\|_2=\|u_g - v_g\|_2$, by the definition of $Tr$ it follows that 
$$
\|\pi_g(e_A)- e_A\|_{2, Tr} = \|u_g - v_g\|_2 \leq \delta
$$
Thus, the element $\xi_0$ 
of minimal norm $\| \ \|_{2,Tr}$ in the weakly compact  convex set obtained by taking weak closure of the convex set $\text{\rm co}\pi_\Gamma(e_A)=\text{\rm co} \{ u_gv_g^* e_A v_gu_g^*\mid g\in \Gamma \}$,  
is fixed by $\pi_\Gamma$ and satisfies $\|\xi_0- e_A \|_{2,Tr}\leq \varepsilon^2/4$. This element lies in $\tilde{A}_+$ with $0\leq \xi_0\leq 1$. 
Its spectral projection $e$ corresponding to the interval $[1/2, 1]$ satisfies $\| e_A  - e\|_{2,Tr}\leq \varepsilon$ and it is fixed by $\pi_\Gamma$.  
Also, since $\tilde{A}$ is commutative and $e, e_A$ are projections, we have $\|e_A- e\|_{1,Tr}=\| e_A  - e\|_{2,Tr}^2\leq \varepsilon^2$.  
Since $\Phi(e_A)=1$, this also shows that $|Tr(e)-1|\leq \varepsilon^2$. 

Thus, $b=\Phi(e) \in M_+$ commutes with $u_g, g\in \Gamma$ and $|\tau(b)-1|\leq \alpha^2$.  Moreover, since $e\in \tilde{A}$ and $\Phi(\tilde{A})=A$, 
it follows that $b=\Phi(e)$ lies in $A$. 
But since $\Gamma$ acts ergodically on $A$, $\{u_g \mid g\in \Gamma\}$ has trivial relative commutant in $A$. 
Hence, $b$ is a scalar $\varepsilon^2$-close to $1$, with $\varepsilon < 1$. But a projection $e\in\tilde{A}$ that has trace $<2$ and expects on 
a scalar in $M$, must have trace $1$ and be of the form $e=ve_Av^*$ for some $v\in \Cal N_M(A)$ with $\|v-1\|_2 = \|ve_Av^*-e_A\|_{2,Tr} =\|e-e_A\|_{2,Tr} \leq \varepsilon$. 

The relation $\pi_g(e)=e$ amounts to $\pi_g(ve_Av^*)=ve_Av^*$, thus $v^*u_gvv_g^*e_Av_gv^*u_g^*v$ $=e_A$ meaning that Ad$(u_g)=\text{\rm Ad}(vv_gv^*)$, $\forall g\in \Gamma$, 
i.e., Ad$(v)$ intertwines $\sigma_g, \theta^{-1}\rho_g \theta$, with $\|\text{\rm Ad}(v)-1\|_2=\|ve_Av^*-e_A\|_{2,Tr}\leq \varepsilon$. 
\hfill $\square$

\proclaim{3.11.  Corollary}  Let $\Gamma$ be a property $(\text{\rm T})$ group. 

$1^\circ$ Two free ergodic pmp $\Gamma$-actions are app-conjugate if and only if they are conjugate.  

$2^\circ$ Let $\Gamma \curvearrowright^\sigma (X,\mu)$ be a pmp $\Gamma$-action. If $(\theta_n)_n \subset \text{\rm Aut}(X,\mu)$ 
are so that $\lim_n \|[\theta_n, \sigma_g]\|_2=0$, $\forall g\in \Gamma$, then there exist $(\theta'_n)_n \subset \text{\rm Aut}(X,\mu)$ 
such that $[\theta'_n, \sigma_g]=0$, $\forall n, g$, and $\lim_n \|\theta'_n - \theta_n\|_2=0$. 
\endproclaim
\noindent
{\it Proof}. $1^\circ$ Let $F\subset \Gamma$ be a finite set of generators. 
Approximate conjugacy of two free ergodic pmp actions $\Gamma \curvearrowright^\sigma X$, $\Gamma \curvearrowright^\rho Y$ 
means there exists a sequence of measure preserving isomorphisms $\theta_n: X  \simeq Y$ such that if we denote 
$\delta_n= \max_{g\in F} \|\theta_n \sigma_g \theta_n^{-1}-\rho_g\|_2$, then $\delta_n \rightarrow 0$.  By Lemma 3.10, there exist measure preserving  
isomorphisms $\theta'_n: X\simeq Y$ such that $\lim_n \|\theta'_n - \theta_n\|_2 =0$ and $\theta'_n \sigma_g = \rho_g \theta'_n$, $\forall g\in \Gamma$. 

$2^\circ$ This is just Lemma 3.10 applied in the case $\rho_g = \sigma_g$ and $\theta=\theta_n\in {\text{\rm Aut}}(X,\mu)$, 
as $n \rightarrow \infty$. 
\hfill $\square$

\vskip .1in

\noindent
{\bf 3.12. Remark}. $1^\circ$ Note that the spectral behavior of $\Gamma \curvearrowright L^2(A^\omega)$ is an $\omega$-conjugacy invariant for 
$\Gamma \curvearrowright (A,\tau)$. In particular, strong ergodicity of $\Gamma \curvearrowright A$  
and weak mixingness of $\Gamma \curvearrowright A^\omega$ are $\omega$-conjugacy invariants. 
Also, the property of having a $\Gamma$-invariant subalgebra $A_0 \subset  A^\omega$ satisfying certain properties, 
like relative property (T) in the sense of (Sec. 4 in [P5]), is of course $\omega$-conjugacy invariant. 

\vskip .05in
$2^\circ$ Another app-conjugacy  invariant for an action 
$\Gamma \curvearrowright^\sigma (A,\tau)$ is its {\it app $\omega$}-{\it centralizer} group  $\text{\rm Aut}^\sigma_{app,\omega}(A)$, defined as the quotient between  
the group of  sequences of automorphisms $(\theta_n)_n$ of $(A,\tau)$ 
that $\omega$-asymptotically  commute with $\sigma_g, \forall g \in \Gamma$, i.e., $\lim_\omega \| [\sigma_g, \theta_n]\|_2$ $=0, \forall g$, 
by the subgroup of sequences of automorphisms $(\theta_n)_n$ 
satisfying $\lim_\omega \|\theta_n - id_A\|_2=0$. This group coincides with the group of automorphisms $\theta$ of $A^\omega$ that are 
of ultraproduct form, $\theta=(\theta_n)_n$, and commute 
with $\sigma^\omega$. 

Note that this group can be identified with a subgroup of the group of automorphisms of $M(\omega)$ that 
leave $A^\omega$ invariant and are of ultraproduct form when restricted to $A^\omega$. As such, it has a topology inherited from 
Aut$(A^\omega \subset M(\omega))$. But we also endow the  group $\text{\rm Aut}^\sigma_{app,\omega}(A)$ with the topology given by the distance 
$\|(\theta_n)_n - (\theta'_n)_n\|_2:=\lim_{n \rightarrow \omega} \|\theta_n - \theta'_n\|_2$. 

Note that  $\text{\rm Aut}^\sigma_{app,\omega}(A)$ contains the group $(\text{\rm Aut}^\sigma(A))^\omega$, of sequences 
of automorphisms $(\theta_n)_n$ of $(A,\tau)$, with all $\theta_n$ in  
the centralizer of $\sigma$ in Aut$(A)$, $\text{\rm Aut}^\sigma(A)$. Moreover, the restriction of the distance on $\text{\rm Aut}^\sigma_{app, \omega}(A)$ 
to this subgroup gives the discrete $\{0, \sqrt{2}\}$-valued distance. 

An interesting result here would be to show that in certain cases 
we have $\text{\rm Aut}^\sigma_{app,\omega}(A)=(\text{\rm Aut}^\sigma(A))^\omega$. 
This amounts to proving that if  $\lim_\omega \|[\sigma_g, \theta_n]\|_2 =0$, $\forall g$, then there exist $\theta'_n \in 
\text{\rm Aut}^\sigma(A)$ such that $\lim_\omega \|\theta_n - \theta'_n\|_2=0$. 
Note that 3.11.2$^\circ$ shows  that for a group $\Gamma$  with the property (T), this is indeed the case.

\vskip .05in
$3^\circ$ We do not know whether $\omega$-conjugacy for property (T) groups is weaker than conjugacy. In fact, $\omega$-conjugacy 
of actions seems to be a rather esoteric, hard to understand concept, for which it seems unlikely that an entropy invariant can be developed 
and whose symmetry groups seem difficult to calculate (see also 5.6.2).

\vskip .05in
$4^\circ$ It is not known whether every non-amenable group $\Gamma$ has at least two non weakly conjugate free ergodic pmp actions. 
But it has been noted in [AE] that strong ergodicity is a weak conjugacy invariant, and thus if $\Gamma$ does not have property (T) 
then by ([Sc2], [CW]) it  does have two non weakly conjugate free ergodic pmp actions (see 4.9 for a strengthening of this result). 
The general belief is that in fact any non-amenable group 
$\Gamma$ has continuously many non weakly conjugate actions. Note that by [AE], this has been checked for many 
groups $\Gamma$, including free groups $\Bbb F_n, 2\leq n \leq \infty$, 
and linear groups with the property (T). Since the app and $\omega$-conjugacies are stronger than weak-conjugacy, each one of these 
examples of non weak conjugate actions of a fixed group $\Gamma$ gives examples of  actions that are not app/$\omega$-conjugate either. 
Since for property (T) groups app-conjugacy is the same as conjugacy, it follows that any non-amenable group has at least two non app-conjugate 
actions. But we do not know how to prove the same for $\omega$-conjugacy, nor how 
to show that any non-amenable group has uncountably many 
non app/$\omega$-conjugate free ergodic pmp actions.

\heading 4.  Approximate orbit equivalence 
\endheading

The description of app-conjugacy of group actions in terms of their ultrapowers, as well as the notion of $\omega$-conjugacy, 
lead to two weak versions of the  orbit equivalence of group actions (respectively of the equivalence of Cartan inclusions) that we discuss in this section. 
We also introduce a notion of weak orbit equivalence of group actions (engendered by weak conjugacy) and relate it with these other concepts. 
\vskip .1in 

\noindent
{\bf 4.1. Definition.} Let $\omega$ be a free ultrafilter on $\Bbb N$. Two pmp actions $\Gamma \curvearrowright^\sigma X$, $\Lambda \curvearrowright^\rho Y$ 
of countable groups $\Gamma, \Lambda$ are 
$\omega$ - {\it orbit equivalent} ($\omega$-{\it OE}), if there exists an isomorphism $\theta$ of Cartan inclusions $(A^\omega \subset M(\omega)) 
\simeq (B^\omega \subset N(\omega))$, where we have denoted $A=L^\infty(X)$, $B=L^\infty(Y)$, 
$M= L(\Cal R_\sigma)$, $N= L(\Cal R_\rho)$, $M(\omega)=A^\omega \vee \Cal N_M(A)\subset M^\omega$, 
$N(\omega)=B^\omega \vee \Cal N_N(B)\subset N^\omega$. 

Note that by Proposition 2.6, if we denote by 
$[\sigma^\omega]$, $[\rho^\omega]$ the full groups of ultrapower actions $\Gamma\curvearrowright^{\sigma^\omega} A^\omega$, 
$\Lambda\curvearrowright^{\rho^\omega} B^\omega$ (or alternatively, of the 
the Cartan inclusions $A^\omega \subset M(\omega)$, 
$B^\omega \subset N(\omega)$), then this condition is equivalent to the existence of an isomorphism 
$\theta: A^ \omega \simeq B^\omega$ that intertwines their full groups, $\theta([\sigma^\omega])\theta^{-1}=[\rho^\omega]$. 

More generally, two Cartan inclusions $A\subset M$, $B\subset N$ are $\omega$-{\it equivalent} 
if there exists an isomorphism $\theta$ of the Cartan inclusions $(A^\omega \subset M(\omega)) 
\simeq (B^\omega \subset N(\omega))$.

\vskip .1in 

Like for $\omega$-conjugacy, we do not know whether 
the $\omega$-OE property depends on the free ultrafilter $\omega$. However, we will next show that the strengthening 
of this condition requiring that 
the isomorphism $\theta: A^\omega \simeq B^\omega$ intertwining the full groups $[\sigma^\omega]$, $[\rho^\omega]$, 
be of ultrapower form $\theta=(\theta_n)_n$, where $\theta_n: A \simeq B$, $\forall n$, does not depend on $\omega$. 
We will do this by showing that this condition is equivalent to a local condition which can be viewed 
as the OE-type equivalence relation for group actions that's entailed by app-conjugacy, and which we will thus call approximate orbit equivalence.

\proclaim{4.2. Proposition} Let $\Gamma \curvearrowright X$, $\Lambda \curvearrowright Y$ be pmp actions of countable groups. 
Denote $L^\infty(X)=A\subset M=L(\Cal R_{\Gamma \curvearrowright X})$, $L^\infty(Y)=B\subset N=L(\Cal R_{\Lambda \curvearrowright Y})$ the  
associated Cartan inclusions. The following properties are equivalent:

$(a)$ There exists a free ultrafilter $\omega$ on $\Bbb N$ for which the following property holds: 
there exists an isomorphism $ \theta:(A^\omega \subset M(\omega)) \simeq (B^\omega \subset N(\omega))$
whose restriction to $A^\omega$ is of ultra product form, $\theta=(\theta_n)_n$, with $\theta_n : A \simeq B$, $\forall n$.   

$(b)$ The property $(a)$ holds true for any free ultrafilter $\omega$ on $\Bbb N$.

$(c)$ There exist maps $\Gamma \ni \phi \mapsto t^\phi = \{t^\phi_\psi\}_{\psi \in \Lambda}$  
and $\Lambda \ni \psi \mapsto s^\psi = \{s^\psi_\phi\}_{\phi \in \Gamma}$, with $\Sigma_\phi s_\phi^\psi = \Sigma_\psi t_\psi^\phi=1$, 
$0\leq s^\psi_\phi, t^\phi_\psi \leq 1$, $\forall \phi, \psi$, 
such that for any $E\subset \Gamma$, $F\subset \Lambda$ finite and any $\varepsilon > 0$, there exist an isomorphism $\theta': A \simeq B$ and 
mutually orthogonal projections $\{p^{\psi_0}_\phi\}_{\phi \in \Gamma}\subset A$, resp. $\{q^{\phi_0}_\psi\}_{\psi \in \Lambda} \subset B$, 
so that $\{\phi(p_\phi^{\psi_0})\}_\phi$ $($resp. $\{\psi(q^{\phi_0}_\psi)\}_\psi$ $)$ are mutually disjoint as well, 
$\tau(p^{\psi_0}_\phi)\leq s^{\psi_0}_\phi$, $\tau(q^{\phi_0}_\psi)\leq t^{\phi_0}_\psi$, 
and $\|{\theta'} \phi_0 {\theta'}^{-1} - \oplus_\psi \psi_{|Bq^{\phi_0}_\psi}\|_2 \leq \varepsilon$, 
$\|{\theta'}^{-1} \psi_0 {\theta'} - \oplus_\phi \phi_{|Ap^{\psi_0}_\phi}\|_2 \leq \varepsilon$, 
$\forall \phi_0\in E, \psi_0\in F$. 

$(d)$ There exist a sequence of isomorphisms $\theta_n : A \simeq B$ and sequences of mutually orthogonal projections $\{p^\psi_{\phi,n}\}_\phi \subset \Cal P(A)$, for $\psi\ \in \Lambda$,
$($respectively  $\{q^\phi_{\psi,n}\}_\psi$, for $\phi\in \Gamma)$ 
with  $\{\phi(p_{\phi,n}^{\psi})\}_\phi$ $($resp. $\{\psi(q^{\phi}_{\psi,n})\}_\psi$ $)$ mutually orthogonal, such that 
$$
\lim_n \|\theta_n \phi \theta_n^{-1} - \oplus_\psi \psi_{|Bq^\phi_{\psi,n}}\|_2=0, \forall \phi\in \Gamma
$$
$$
\lim_n \|\theta_n^{-1} \psi \theta_n - \oplus_\phi \phi_{|Ap^\psi_{\phi,n}}\|_2=0, \forall \psi\in \Lambda
$$
$$
\Sigma_\phi \liminf_n \tau(p^\psi_{\phi,n})=1, \Sigma_\psi \liminf_n \tau(q^\phi_{\psi,n})=1, \forall  \psi \in \Lambda, \phi\in \Gamma.
$$ 

\endproclaim
\noindent
{\it Proof}. Let $\{u_\phi \mid \phi \in \Gamma \} \subset M=L([\Gamma])=L(\Cal R_\Gamma)$ 
(respectively $\{v_\psi \mid \psi \in \Lambda\}\subset N=L([\Lambda])=L(\Cal R_\Lambda)$)  
denote the canonical unitaries in $M$ (resp. $N$), implementing the action $\Gamma \curvearrowright A$ (resp. $\Lambda \curvearrowright B$). 
 
Assume $(a)$ holds true. Then for each  $\phi_0\in \Gamma, \psi_0\in \Lambda$, there exist partitions of $1$ with projections 
$\{Q^{\phi_0}_\psi\}_\psi\subset \Cal P(B^\omega)$, $\{P^{\psi_0}_\phi\}_\phi\subset \Cal P(A^\omega)$, 
such that 

$$\theta(u_{\phi_0})=(\Sigma_\psi v_\psi Q^{\phi_0}_\psi)b_{\phi_0}, 
\theta^{-1}(v_{\psi_0})=(\Sigma_\phi u_\phi P^{\psi_0}_\phi)a_{\psi_0},\tag 4.2.1
$$ 
for some $a_{\psi_0}\in \Cal U(A^\omega)$, $b_{\phi_0}\in \Cal U(B^\omega)$. 

Denote $s^{\psi_0}_\phi=\tau(P^{\psi_0}_\phi)$, $t^{\phi_0}_\psi=\tau(Q^{\phi_0}_\psi)$, then choose projections 
$P^{\psi_0}_{\phi, m}\in \Cal P(A)$, $Q^{\phi_0}_{\psi, m}\in \Cal P(B)$ so that $P^{\psi_0}_\phi=(P^{\psi_0}_{\phi, m})_m$ 
and $Q^{\phi_0}_\psi=(Q^{\phi_0}_{\psi, m})_m$. Moreover, by taking into account the properties that $P^{\psi_0}_\phi$, 
$Q^{\phi_0}_\psi$ satisfy, it follows that we can make this choice such that for each $m$ we have  
$\tau(P^{\psi_0}_{\phi,m})\leq s^{\psi_0}_\phi$, $\tau(Q^{\phi_0}_{\psi, m})\leq t^{\phi_0}_\psi$ and $\Sigma_\phi P^{\psi_0}_{\phi, m} \leq 1$, 
$\Sigma_\phi u_{\phi}P^{\psi_0}_{\phi, m}u_{\phi}^* \leq 1$, $\Sigma_\psi Q^{\phi_0}_{\psi,m} \leq 1$, $\Sigma_\psi v_\psi(Q^{\phi_0}_{\psi,m})v_\psi^* \leq 1$. 

But then, if $E\subset \Gamma$, $F\subset \Lambda$ are finite sets and $\varepsilon >0$, there exists $m$ ``close to $\omega$'' 
so that $p^{\psi_0}_\phi=P^{\psi_0}_{\phi,m}$ and $q^{\phi_0}_{\psi}=Q^{\phi_0}_{\psi,m}$ satisfy all the conditions in $(c)$. 

To see that $(c)$ implies $(d)$, choose sequences of finite subsets $E_n \nearrow \Gamma$, $F_n \nearrow \Lambda$ 
and for each $n$ apply $(c)$ to $E_n, F_n$, $\varepsilon_n=2^{-n}$ to get $\{p^{\psi_0}_{\phi, n}\}_\phi\subset \Cal P(A)$, $\{q^{\phi_0}_{\psi, n}\}_\psi\subset \Cal P(B)$ 
and $\theta_n : A \simeq B$ satisfying the conditions in $(c)$ for all $\phi_0\in E_n$, $\psi_0\in F_n$. Then $\{\theta_n\}_n$, $\{p^\psi_{\phi,n}\}_n$, $\{q^\phi_{\psi,n}\}_n$ 
clearly satisfy $(d)$.

Finally, assume $(d)$ holds true and let $\omega$ be a given (but arbitrary) free 
ultrafilter on $\Bbb N$.  Denote $\theta'=(\theta_n)_n: A^\omega \simeq B^\omega$. Let 
$P^{\psi_0}_\phi=(p^{\psi_0}_{\phi, n})_n  
\in \Cal P(A^\omega)$, $Q^{\phi_0}_\psi=(q^{\phi_0}_{\psi, n})_n\in \Cal P(B^\omega)$ and note that 
by the last condition in $(d)$ we have $\Sigma_\phi \tau(P^{\psi_0}_\phi) = \Sigma_\phi \lim_\omega \tau(p^{\psi_0}_{\phi,n})=1$,  
and thus $\Sigma_\phi P^{\psi_0}_\phi =1$, $\forall \psi_0$. Similarly,  $\Sigma_\psi Q^{\phi_0}_\psi=1$, $\forall \phi_0$. 
Moreover, by conditions at the beginning of $(d)$ we also get $\Sigma_\phi \phi(P^{\psi_0}_\phi)=1$, $\Sigma_\psi \psi(Q^{\phi_0}_\psi)=1$. 

If we now let $\theta(u_{\phi_0})=\Sigma_\psi v_\psi Q^{\phi_0}_\psi$, $\forall \phi_0\in \Gamma$, then it is easy to check that there exists a 
unique isomorphism $\theta: M(\omega)\simeq N(\omega)$ which on $A^\omega$ acts as $\theta'$ and 
on $\{u_\phi\}_\phi$ takes these assigned values. This shows that $(d)$ implies $(b)$ and thus finishes the proof. 

\hfill $\square$

\vskip .1in 

\noindent
{\bf 4.3. Definition.} Two pmp actions of countable groups $\Gamma \curvearrowright X$, $\Lambda \curvearrowright Y$ 
are  {\it approximately orbit equivalent} ({\it app-OE})  if any of the equivalent conditions in 4.2 is satisfied.

\vskip .1in 

\noindent
{\bf 4.4. Remark.} We  define {\it app-equivalence} of separable Cartan inclusions $A_i\subset M_i$, $i=1,2$, 
by viewing them as pairs $(\Cal G_i, v_i/{\sim})$,  
(a full group and a 2-cocycle, cf. Proposition 2.6), then noticing that analogues of $(a)$, $(b)$, $(d)$ in 4.2 above are equivalent, where $(a)$, $(b)$ are the same, 
and $(d)$ adds the condition that $\theta_n$ asymptotically intertwine the 2-cocycles $v_1/{\sim}, v_2/{\sim}$, and asking that any of these 
equivalent conditions holds true. If $v_1=1=v_2$, then the Cartan inclusions arise  
from  pmp actions of countable groups $\Gamma_i \curvearrowright A_i$, and the app-equivalence of $A_i\subset M_i$, $i=1,2$, amounts to app-OE  
of these group actions, in the sense of 4.3.

\proclaim{4.5. Lemma} Let $A\subset M$, $B\subset N$ be Cartan inclusions with $M, N$ separable. Let $\Cal U\subset \Cal N_M(A)$  
be a countable $\| \ \|_2$-dense subgroup, with $\Cal U_0=\Cal U\cap A$ normal in $\Cal U$ and dense in $\Cal U(A)$, and let $v_0=1, v_1, ... \in q\Cal N_N(B)$ 
be an orthonormal basis of $N$ over $B$. The following conditions are equivalent: 

$(a)$ There exists a non-degenerate Cartan embedding of $A\subset M$ into $B^\omega\subset N(\omega)$ for some free ultrafilter $\omega$ on $\Bbb N$. 

\vskip .05in 

$(b)$  There exists a non-degenerate Cartan embedding of $A\subset M$ into $B^\omega 
\subset N(\omega)$ for any free ultrafilter $\omega$ on $\Bbb N$.

\vskip .05in

$(c)$ The following property, which we denote $(A\subset M) \prec_w (B\subset N)$, holds true: 

\vskip .05in 
\noindent
There exists a map $\Cal U \ni u \mapsto t^u = \{t^u_l\}_{l\geq 0}$,   
with $\Sigma_l t_l^u=1$, $0\leq  t^u_l\leq 1$, $\forall u, l$, such that 
for any $L\geq 1$, any $u_1, ..., u_n \in \Cal U$, any $\varepsilon > 0$ and $K\geq 1$, there exist  $u_{n+1}, ..., u_m\in \Cal U$, 
$\{q^i_l\}_{l\geq 0}\subset \Cal P(B)$, with $\Sigma_l q^i_l, \Sigma_l v_lq^i_l v_l^* \leq 1$, $\tau(q^{i}_l)\leq t^{u_i}_l$, $\forall l\geq 0, 1\leq i \leq m$, $q^i_0=1$ 
whenever $u_i\in \Cal U_0$, 
and $b_1, ..., b_m \in \Cal U(B)$ such that  if we denote $u'_i=(\Sigma_l v_l q_l^i)b_i$, $1\leq i \leq m$, then 
$v_l \in_\varepsilon \Sigma_{i=1}^m u'_i B$, $\forall 0\leq l \leq L$, and given any word $w$ of length $\leq K$ in $m$ letters $($and their formal 
inverses$)$, we have 
$|\tau (w(\{u'_i\}_i))-\tau(w(\{u_i\}_i))| \leq \varepsilon$. 
\endproclaim
\noindent
{\it Proof}.  $(a)\implies (c)$ Denote by $\theta$ the isomorphism of $M$ into $N(\omega)$ taking $A$ into $B^\omega$ and $\Cal N_M(A)$ 
into $\Cal N_{N(\omega)}(B^\omega)$. 

Thus, for each $u\in \Cal U$ there exist $\{q^u_l\}_{l\geq 0} \subset \Cal P(B^\omega)$, $b^u\in \Cal U(B^\omega)$ 
such that: 

$$ \Sigma_l q^u_l = 1 = \Sigma_l v_lq^u_lv_l^*; \tag 1
$$
$$
 \theta(u)=(\Sigma_l v_l q^u_l)b^u, \forall u\in \Cal U;  \tag 2
 $$
 $$
\{v_l\}_l \subset \overline{\text{\rm sp}\theta(\Cal U)B^\omega}. \tag 3
$$ 

Denote $t^u_l=\tau(q^u_l)$ and represent each $q^u_l$ as $(q^u_{l,k})_k$ with $q^u_{l,k}\in \Cal P(B)$, $\tau(q^u_{l,k})\leq t^u_l$, 
$\Sigma_l q^u_{l,k} \leq 1$, $\Sigma_l v_lq^u_{l,k}v_l^* \leq 1$, $\forall u\in \Cal U$, $\forall k\geq 1$. Also, let $b^u=(b^u_k)_k$ with $b^u_k\in \Cal U(B)$.  

Take now $u_1, ..., u_n \in \Cal U$, $L\geq 1$ and $\varepsilon >0$. Note first that by (3), there exists $u_{n+1}, ..., u_m\in \Cal U$ 
such that $v_l\in_{\varepsilon/2} \Sigma_{i=1}^m \theta(u_i) B^\omega$, $\forall l\leq L$. Thus, if $K>>1$ is also given, 
then there exists some $k\geq 1$ ``close to $\omega$'',  such that if for each $1\leq i \leq m$ we denote $q^i_l=q^{u_i}_{l,k}$, $b_i=b^{u_i}_k$, 
then all conditions in $4.6.(c)$ are satisfied. 

To see that $(c)$ implies $(b)$, choose first an enumeration $\Cal U=\{u_j\}_{j \geq 1}$ and for each $n$ apply 
$4.6.(c)$ to $u_1, ..., u_n$,  $F=\{v_0, v_1, ..., v_n\}$ and $\varepsilon=2^{-n}$, to get  some larger integer $k_n\geq n$, projections 
$\{q^j_{l,n}\in \Cal P(B) \mid 1\leq j \leq k_n, l\geq 0\}$, and $b_{j,n}\in \Cal U(B)$, $1\leq j \leq k_n$, such that $q^j_0=1$ whenever $u_j \in \Cal U_0$, 
$\Sigma_{l\geq 0} q^j_{l,n}, \Sigma_l v_lq^j_{l,n} v_l^* \leq 1$, $\tau(q^{j}_{l,n})\leq t^{u_j}_l$ and such that  if we denote 
$u'_{j,n} = (\Sigma_l v_lq^j_{l,n})b_{j,n}$ then 
any word $w$ of length $\leq 2^n$ in $k_n$ letters  satisfies the moments condition $|\tau(w(\{u'_{j,n}\}_{1\leq j \leq k_n}))-\tau(w(\{u_j\}_{1\leq j \leq k_n}))|\leq 2^{-n}$. 

Denote $q^j_l=(q^j_{l,n})_n \in \Cal P(B^\omega)$, $b_j=(b_{j,n})_n \in \Cal U(B^\omega)$, $j\geq 1$, $u'_j=(u'_{j,n})_n = (\Sigma_{l \geq 0} v_l q^j_l)b_j \in N^\omega$. 
Note that, by the definitions, $u'_j$ lie in the closure of $\Sigma_{l \geq 0} v_l B^\omega$, i.e. in $N(\omega)$. Moreover,  
the moments condition implies that the map $\theta: \Cal U \rightarrow N(\omega)$ defined by $\theta(u_j)=u'_j, j\geq 1,$ extends to  a trace-preserving 
isomorphism of $M$ into $N(\omega)$ with $\theta(\Cal N_M(A))$ $\subset \Cal N_{N(\omega)}(B^\omega)$. Also, the condition $q^j_0=1$ implies that 
$\theta(\Cal U_0)\subset B^\omega$. Thus, $\theta$ implements a commuting square embedding of $A\subset M$ into $B^\omega \subset N(\omega)$, 
which by the condition $v_l \in \overline{\Sigma_i u'_i B^\omega}$, $\forall l\geq 0$, follows non-degenerate. 

\hfill $\square$

\proclaim{4.6. Corollary} Let $A\subset M$, $B\subset N$ be Cartan inclusions of separable von Neumann algebras. With 
the same notations as in $4.2$, the following conditions are equivalent: 

$(a)$ There exists a free ultrafilter $\omega$ on $\Bbb N$ for which 
there exist non-degenerate Cartan embeddings of $A\subset M$ into $B^\omega \subset N(\omega)$ and of $B\subset N$ into $A^\omega\subset M(\omega)$. 

$(b)$ Property $(a)$ holds true for any free ultrafilter $\omega$ on $\Bbb N$. 

$(c)$ With the notation in $4.5$, we have $(A\subset M)\prec_w (B\subset N)$ and $(B\subset N)\prec_w (A\subset M)$.

\endproclaim
\noindent
{\it Proof}. This is now trivial by 4.5. 
\hfill $\square$

\vskip .1in
\noindent
{\bf 4.7. Definition.} Two pmp actions of countable groups $\Gamma \curvearrowright X$, $\Lambda \curvearrowright Y$ 
with their associated Cartan inclusions $A\subset M=L(\Cal R_{\Gamma \curvearrowright X})$, 
$B \subset N=L(\Cal R_{\Lambda \curvearrowright Y})$ (respectively two arbitrary separable 
Cartan inclusions $A\subset M$, $B\subset N$) are {\it weakly orbit equivalent}, abbreviated {\it weakly-OE} (respectively {\it weakly equivalent})  if 
any of the equivalent conditions in 4.6 above holds true.

\vskip .1in

It is clear that OE $\implies$ app-OE 
$\implies$ $\omega$-OE $\implies$ weak-OE and that app-OE is implied by app-conjugacy, $\omega$-OE is implied by $\omega$-conjugacy and 
weak-OE is implied by weak-conjugacy. Clearly $\omega$-OE and app-OE are equivalence relations. The characterization $4.5.(c)$ 
is easily seen to imply that the subordination relation $\prec_w$ is transitive, implying that weak-OE (respectively weak-equivalence) is an equivalence relation 
for group actions (resp. for Cartan inclusions).

\proclaim{4.8.  Proposition} $1^\circ$ Strong ergodicity is a weak-OE invariant $($thus an app-OE and $\omega$-OE invariant as well$)$. 

$2^\circ$ Amenability, relative Haagerup property, co-rigidty are weak-OE invariants $($thus app-OE and $\omega$-OE invariants as well$)$. 

$3^\circ$ Cost is a weak-OE invariant $($thus  app-OE and $\omega$-OE invariant as well$)$. 

$4^\circ$ $L^2$-Betti numbers  are weak-OE invariants $($thus app-OE and $\omega$-OE invariants as well$)$. 
\endproclaim
\noindent
{\it Proof}. Parts $1^\circ$ and $3^\circ$ are trivial by the characterization 4.6$(a)$ of weak-OE and Propositions 2.10 respectively 2.12.$2^\circ$. 

The weak-OE invariance of Haagerup property and of co-rigidity in part $2^\circ$ are consequences of 4.6$(a)$ and of Proposition 2.13, while the weak-OE invariance 
of amenability is trivial, for instance by using [CFW]. 

Part $4^\circ$ is an immediate consequence of 4.6$(a)$ and  2.12.3$^\circ$. 
\hfill $\square$

\proclaim{4.9.  Corollary} If $\Gamma$ is non-amenable and does not have property $\text{\rm (T)}$, then $\Gamma$ 
has two non weak-OE $($thus also non $\omega$-OE and non app-OE$)$ free ergodic pmp actions.  
\endproclaim
\noindent
{\it Proof}. By [Sc2] any Bernoulli action of a non-amenable group is strongly ergodic, 
while by [CW], if $\Gamma$ does not have property (T) then it has a free ergodic pmp action that is not strongly ergodic, and so the statement 
follows from 4.8.$1^\circ$.  
\hfill $\square$

\vskip .1in 
\noindent
{\bf 4.10. Remarks.} $1^\circ$ Note that the property $(A\subset M) \prec_w (B\subset N)$ in  4.5$(c)$ requires much more than just the fact that ``$A\subset M$ can be simulated 
inside $B\subset N$''. Indeed, such a {\it simulation condition} would only require  
\vskip .02in
{\it $\forall v_1, ..., v_n \in q\Cal N_M(A)$, $\varepsilon > 0$, $K\geq 1$,  
$\exists v'_1, ..., v'_n \in q\Cal N_N(B)$, such that given any word $w$ of length $\leq K$ in $n$ letters, we have $|\tau (w(\{v_i'\}_i))-\tau(w(\{v_i\}_i))| \leq \varepsilon$. }
\vskip .02in
In addition to this, the definition of $\prec_w$ in 4.5$(c)$ requires that: $(1)$ the ``simulating elements'' $v_i'$ are obtained from the initial group elements 
by patching them with uniform weights (independent of $\varepsilon$); $(2)$ the simulating elements $v'_i$ tend to exhaust  $N$ (as $\varepsilon \rightarrow 0$).  We have encountered 
the ``uniform weights''  condition also in the definition of app-OE and it will appear again when defining app-cocycles in Sec. 5.6.  This corresponds to 
the fact that $v'_i$ come from the normalizer of $B^\omega$ in $N(\omega)$ (not from the normalizer of $B^\omega$ in $N^\omega$). 

The above simulation condition, which we will denote $(A\subset M)\prec_{CAE} (B\subset N)$,  
is a version for Cartan inclusions of  ``Connes approximate embedding'' (abbreviated {\it CAE}) subordination between II$_1$ factors, $M\prec_{CAE} N$, 
which requires that  $M$ can be ``simulated'' inside $N$, in moments, with the corresponding CAE equivalence $M\sim_{CAE} N$ requiring $N\prec_{CAE} M$ as well. 
If $D\subset R$ 
denotes the (unique) Cartan MASA in the hyperfinite II$_1$ factor, then   one has a 
(plain) Cartan embedding of $(D\subset R)$ into any Cartan inclusion $(B\subset N)$ 
with $N$ a II$_1$ factor, so in particular $(D\subset R) \prec_{CAE} (B\subset N)$. 
The question of whether $(B\subset N) \prec_{CAE} (D\subset R)$ for a Cartan inclusion coming from a free ergodic pmp action 
of a countable group $\Lambda$, $B\subset N=B\rtimes \Lambda$, amounts to asking whether the {\it action is sofic}. While asking 
this to be true for some free $\Lambda$-action (equivalently for a Bernoulli $\Lambda$-action, cf. [EL], or 5.1 in [P11]), 
amounts to $\Lambda$ being  a {\it sofic group}. 

If one denotes by $\Cal N(\omega)$ 
the normalizer of $\text{\bf D}(\omega)=\Pi_\omega D_n$ in $\text{\bf M}(\omega)=\Pi_\omega M_{n \times n}(\Bbb C)$,  
where $D_n$ denotes the diagonal subalgebra in $M_{n \times n}(\Bbb C)$, 
then this is easily seen to be equivalent to the existence 
of an embedding $\Lambda \subset \Cal N(\omega)$ with $\Lambda \curvearrowright \text{\bf D}(\omega)$ free, 
for some (equivalently any) free ultrafilter $\omega$ on $\Bbb N$ . 
Such an embedding of a sofic group into  $\Cal N(\omega)$ is called a  
{\it sofic approximation} of $\Lambda$. 

Recall that there are no known examples of groups, 
or actions,  that are not sofic. 
So apriori, it may be that all Cartan inclusions $B\subset N$ are CAE equivalent to $D\subset R$, the same way 
we may have $M\sim_{CAE} R$, for any II$_1$ factor $M$ (cf.  {\it CAE conjecture}). 

\vskip .05in 
$2^\circ$  One can define yet another subordination property between Cartan inclusions $A\subset M$ and $B\subset N$, 
weaker than $\prec_w$ in 4.5, but stronger than $\prec_{CAE}$ in $1^\circ$ above, by requiring 
the existence of a Cartan embedding of $(A\subset M)$ into $(B^\omega \subset N(\omega))$, in the sense of 1.1 in [P11] (so 
an embedding of $M$ into $N(\omega)$, with $A$ taken into $B^\omega$, $\Cal N_M(A)$ into $\Cal N_{N(\omega)}(B^\omega)$, but without the  
non-degeneracy condition sp$MB^\omega$ dense in $N(\omega)$). A description in local terms of this property  in the spirit of $4.5.(c)$, is as follows 
(where $\Cal U, \{v_l\}_{l \geq 0}$ are as in 4.5):  

\vskip .02in

{\it There exists a map $\Cal U \ni u \mapsto t^u = \{t^u_l\}_{l \geq 0}$,   
with $\Sigma_l t_l^u=1$, $0\leq  t^u_l\leq 1$, $\forall u\in \Cal U, l \geq 0$, such that 
for any $u_1, ..., u_n \in \Cal U$, any $\varepsilon > 0$ and $K\geq 1$, there exist  
$\{q^{i}_l\}_l \subset \Cal P(B)$, with $\Sigma_l q^i_l\leq 1$, $\Sigma_l v_lq^i_l v_l^* \leq 1$, $\tau(q^{i}_l)=t^{u_i}_l$, $\forall l \geq 0, 1\leq i \leq n$, 
such that if we denote $u'_i=\Sigma_l v_lq_l^{i}$, $1\leq i \leq n$, then given any word $w$ of length $\leq K$ in $n$ letters, we have 
$|\tau (w(\{u'_i\}_i))-\tau(w(\{u_i\}_i))| \leq \varepsilon$. }
\vskip .02in

We will denote by $\prec_{ww}$ this subordination relation. Note that hereditarity under $\prec_{ww}$  does hold 
for co-amenability and relative Haagerup property, 
but  not for strong ergodicity (nor for cost, $L^2$-invariants, or rigidity). Moreover, by [AW] and [GL], 
one has the following characterization of amenability in this framework, in the spirit of ``von Neumann's problem'': 
if $\Gamma$ is a countable discrete group and $\Gamma \curvearrowright (X,\mu)$ is 
an arbitrary free ergodic pmp 
$\Gamma$-action, then $\Gamma$ is non-amenable  if and only if $(\Bbb F_2 \curvearrowright [0,1]^{\Bbb F_2})\prec_{ww} 
(\Gamma \curvearrowright X)$, where the $\prec_{ww}$ subordination 
for free ergodic pmp actions means that the corresponding Cartan subalgebras satisfy the $\prec_{ww}$-subordination. 
It would be interesting to provide a direct proof of this result, that would not use [GL] and its percolation methods. Indeed, 
by [E] (see also [BHI]), we see that  such a result is sufficient for deriving existence of many non-OE actions of a given non-amenable 
group $\Gamma$, by using the co-induction technique in [E]. 

The $\prec_{ww}$-subordination 
gives rise to the equivalence relation $(A\subset M) \sim_{ww} (B\subset N)$ between Cartan inclusions $A\subset M, B\subset N$ 
(and free ergodic pmp actions that entail them), 
by requiring that $(A\subset M)\prec_{ww} (B\subset N)$ and $(B\subset N)\prec_{ww} (A\subset M)$, and which 
may be interesting to study in its own right.

\vskip .05in 
$3^\circ$  Yet another notion of ``weak subordination'' between free pmp actions  
$\Gamma \curvearrowright^\sigma X$, $\Gamma \curvearrowright^\rho Y$ has been considered 
in [B2]. It requires that  $\sigma$ can be simulated better and better 
by $\Gamma$-actions that are orbit equivalent to $\rho$.  More in the spirit of this paper, 
this amounts to the existence of a sequence of free pmp $\Gamma$-actions $\Gamma \curvearrowright^{\rho_n} Y_n$ that are OE to $\rho$ 
and such that the ultraproduct action $(\rho_n)_n$ of $\Gamma$ on $\tilde{B}=\Pi_{n \rightarrow \omega} L^\infty(Y_n)$, contains a 
$\Gamma$-invariant subalgebra $B\subset \tilde{B}$ with $\Gamma \curvearrowright B$ isomorphic to $\sigma$.  It is shown in [B2] that 
any two free pmp actions $\sigma, \rho$ of a free group $\Gamma = \Bbb F_n$, $1\leq n \leq \infty$, are subordinated one to the other, in this sense.

\heading 5.  $\omega$-OE rigidity and open problems  
\endheading

In this section we discuss the app-OE and $\omega$-OE versions of the various types of OE-rigidity paradigms for group actions:    
OE strong rigidity, OE-superrigidity, cocycle superrigidity, and calculation of OE symmetry groups,  i.e., the 
automorphism group,  cohomology groups, the fundamental group.

Recall that a free ergodic pmp action $\Gamma \curvearrowright X$ is called OE-superrigid if any OE between $\Gamma \curvearrowright X$ 
and another free pmp action $\Lambda \curvearrowright Y$ comes from a conjugacy.  An OE strong rigidity type result 
is a weaker version of this, deriving automatic conjugacy from the orbit equivalence between $\Gamma \curvearrowright X$ and an action  
$\Lambda  \curvearrowright Y$ belonging to a special class of group actions. 

Instead, we will assume here the app-OE or $\omega$-OE  between a certain group action $\Gamma \curvearrowright X$ and 
another action  $\Lambda \curvearrowright Y$ (arbitrary, or from a specific  class), trying to prove that this automatically entails 
a stronger equivalence between them. The natural stronger equivalence to seek for in this case is app-conjugacy and respectively 
$\omega$-conjugacy. We will provide  below large classes of group actions 
$\Gamma \curvearrowright X$ for which, indeed,  any $\omega$-OE (resp. app-OE) with another free pmp action 
comes from an $\omega$-conjugacy (resp. app-conjugacy). 

We will do this by first proving that such {\it app-OE/$\omega$-OE superrigidity} for a certain free ergodic pmp action $\Gamma \curvearrowright X$ of a group $\Gamma$ 
follows from the fact that all free ergodic pmp $\Gamma$-actions that satisfy certain properties are OE-superrigid. Combining this 
with the striking OE-rigidity results in ([K1, K2], [CK], [MS]),   
we can then derive many examples of such phenomena.

\proclaim{5.1.  Proposition} Assume the countable group $\Gamma$ has the property 
that any free, strongly ergodic pmp action $\Gamma \curvearrowright^\sigma X$ is  OE-superrigid. Then any such action is  
$\omega$-OE superrigid $($resp. app-OE superrigid$)$ as well, i.e., 
any $\omega$-OE $($resp. app-OE$)$ of this action with another free pmp action  $\Lambda \curvearrowright^\rho Y$ comes from an $\omega$-conjugacy 
$($resp. app-conjugacy$)$ of $\Gamma \curvearrowright X$, $\Lambda \curvearrowright Y$. More specifically, if $\theta:A^\omega \simeq B^\omega$ 
implements an $\omega$-OE $($respectively app-OE$)$ between the two group actions, then there exists $u\in \Cal N_{N(\omega)}(B^\omega)$ 
such that $\text{\rm Ad}(u)(\theta(\sigma^\omega(\Gamma))\theta^{-1}) = \rho^\omega(\Lambda)$. 
\endproclaim
\noindent
{\it Proof}. Assume $\Gamma\curvearrowright^\sigma X$ and $\Lambda \curvearrowright^\rho Y$ are $\omega$-OE. 
Thus, with the usual notations $A=L^\infty(X)\subset L^\infty(X)\rtimes \Gamma=M$, 
$B=L^\infty(Y) \subset L^\infty(Y)\rtimes \Lambda=N$, $M(\omega)=A^\omega \rtimes \Gamma$, $N(\omega)=B^\omega \rtimes \Lambda$, 
we have an isomorphism $\theta: (A^\omega \subset M(\omega)) \simeq (B^\omega \subset N(\omega))$. Let $\{u_g \mid g \in \Gamma\}$ 
(respectively $\{v_h \mid h \in \Lambda\} \subset N(\omega)$) 
be the canonical unitaries implementing $\Gamma \curvearrowright A^\omega$ (resp. 
$\Lambda \curvearrowright B^\omega$). Denote by $\Cal Y_0$ 
the set of projections $\{q^g_h \mid g \in \Gamma, h\in \Lambda\}\subset B^\omega$ such that 
$\theta(u_g)=\Sigma_h v_h q^g_h$, $\forall g\in \Gamma$. 

Define now $A_0=A$, $B_0=\vee_{h\in \Lambda} \rho^\omega_h (B \cup \Cal Y_0 \cup \theta(A))$ and for each 
$n\geq 1$ define recursively $A_n=\vee_{g\in \Gamma} \sigma_g^{\omega}(A_{n-1} \cup \theta^{-1}(B_{n-1}))$, $B_n=\vee_{h\in \Lambda} \rho^\omega_h (B_{n-1} \cup \theta(A_n))$. If we then let  
$\tilde{A}=\overline{\cup_n A_n}$, $\tilde{B}=\overline{\cup_n B_n}$,  then  we clearly have $\theta(\tilde{A})=\tilde{B}$, 
$\theta(\{u_g\}_g) \subset \tilde{B}\vee \{v_h\}_h$, 
with $\Gamma$ (resp. $\Lambda$) leaving invariant $\tilde{A}$ (resp. $\tilde{B}$) and acting freely on it. Thus, if we denote $\tilde{M}=\tilde{A}\vee \{u_g\}_g\simeq \tilde{A} \rtimes \Gamma$, $\tilde{N}=\tilde{B}\vee \{v_h\}_h \simeq \tilde{B} \rtimes \Lambda$, 
then the restriction of $\theta$ to $\tilde{M}$ implements an isomorphism of the Cartan inclusions $(\tilde{A}\subset \tilde{M})\simeq (\tilde{B}\subset \tilde{N})$. 

By Proposition 2.10, it follows that $\Gamma \curvearrowright \tilde{A}$ is strongly ergodic (and also free, with $\tilde{A}$ separable). 
The hypothesis implies that $\Gamma \curvearrowright \tilde{A}$ is OE-superrigid. Thus, there exists a unitary element $u\in \Cal N_{\tilde{N}}(\tilde{B})$ 
such that $u\theta(\{u_g \mid g\in \Gamma\})u^* 
\subset \{v_hv \mid h\in \Lambda, v\in \Cal U(B) \}$, where $u_g\in \tilde{M}\subset M(\omega)$, 
$v_h \in \tilde{N} \subset N(\omega)$ are the canonical unitaries. Since $\Cal N_{\tilde{N}}(\tilde{B})\subset \Cal N_{N(\omega)}(B^\omega)$, we have this way shown 
that $\theta$ comes from an $\omega$-conjugacy. Moreover, since $u\in \Cal N_{N(\omega)}(B^\omega)$ can be represented as $(u_n)_n$ 
with $u_n \in \Cal N_N(B)$, it follows that Ad$u=(\text{\rm Ad}u_n)_n$ is of ultraproduct form as well. Thus, if $\theta=(\theta_n)_n$ implements an app-OE, 
then Ad$u  \circ \theta$ implements an app-conjugacy of $\Gamma \curvearrowright X$, $\Lambda \curvearrowright Y$. 
\hfill $\square$

\vskip .1in 

Strictly speaking, we will not be able to 
apply Proposition 5.1 the way it is stated. But its proof 
provides  the template of  how to deduce app-OE/$\omega$-OE superrigidity from OE-superrigidity 
for many ``special'' 
free strongly ergodic 
pmp actions $\Gamma \curvearrowright^\sigma X$. Thus, if say we know that OE-superrigidity holds for all actions of $\Gamma$ that satisfy 
a certain property $\Cal P$, then we can deduce from it app-OE superrigidity for all $\Gamma$-actions $\sigma$ with the property 
that any ``quotient of its ultrapower $\sigma^\omega$'' still satisfies property $\Cal P$.  
This will force us to assume that $\sigma$ satisfies some strengthened form of $\Cal P$, which 
will insure it survives to all (separable) quotients of $\sigma^\omega$ (i.e., $\sigma$ has good $\omega$-{\it permanence}).  

\proclaim{5.2.  Corollary} $1^\circ$  If  $\Gamma=SL(3,\Bbb Z)*_\Sigma SL(3,\Bbb Z)$, where $\Sigma\subset SL(3,\Bbb Z)$ is the subgroup of  matrices $(t_{ij})$ with $t_{31}=t_{32}=0$, 
then any free, strongly ergodic, aperiodic pmp $\Gamma$-action is app-OE superrigid and $\omega$-superrigid.  
\vskip .05in 

$2^\circ$ If $\Gamma$ is of the form $\Gamma = \Gamma_1 *_H \Gamma_2$, 
with $\Gamma_1, \Gamma_2$ lattices in non-compact connected simple Lie groups with trivial center and real rank $\geq 2$, and with 
the common subgroup $H$ non-amenable and satisfying 
$|\Gamma_i/H|=\infty$, as well as the singularity condition $g\in \Gamma_i,  [H: g H g^{-1}\cap H]<\infty \implies g\in H$, 
$i=1, 2$, then any $\Gamma$-action whose restriction to $H$ is both strongly ergodic and aperiodic 
$($e.g., any quotient of a Bernoulli $\Gamma$-action$)$ is app-OE superrigid and $\omega$-superrigid.

\vskip .05in

$3^\circ$ The following $\Gamma$ have the property that any free, aperiodic, strongly  ergodic, pmp 
$\Gamma$-action $($e.g., any quotient of a Bernoulli $\Gamma$-action$)$ is app-OE superrigid and $\omega$-superrigid:   

$(a)$ $\Gamma$ is a finite index subgroup of a finite product of mapping class groups $\Gamma_1, ...,  \Gamma_n$,  
of compact orientable surfaces,   with the  
genus $g_i$ and number of boundary components $p_i$ of $\Gamma_i$ satisfying $3g_i + p_i-4 >0$ and $(g_i, p_i) \neq (2,0), (1,2)$, $\forall i$; 

$(b)$ $\Gamma$ is the pure braid group with $k\geq 2$ strands on a closed, orientable surface of genus $g\geq 2$.

\endproclaim
\noindent
{\it Proof}. Let us first notice that if a free ergodic pmp action $G\curvearrowright^\rho Z$ of a countable group $G$ is 
strongly ergodic and aperiodic, then the restriction of $\rho^\omega: G \rightarrow \text{\rm Aut}(L^\infty(Z)^\omega)$ to any  
$G$-invariant von Neumann subalgebra $C \subset L^\infty(Z)^\omega$ on which it acts freely is still strongly ergodic and aperiodic 
(i.e.,  given any subgroup of finite index $G_0\subset G$, the action $G_0 \curvearrowright C$ is free and strongly ergodic). 
Indeed, this can be trivially seen by the proof of ($1.11(i)$  in [PP]; see also 3.1 in [AE]). 

Now, to prove $1^\circ$, recall first that 
by  (Theorem 1.5(b) in [K2]), any free ergodic pmp action of the group 
$\Gamma=SL(3,\Bbb Z)*_\Sigma SL(3,\Bbb Z)$ is stably (or virtually) OE-superrigid. But the group $\Gamma$ has no finite normal subgroups 
and also, since our $\Gamma$-action is strongly ergodic and aperiodic, by the above remark its restriction 
to any finite index subgroup of $\Gamma$ is still strongly ergodic, so in particular it is ergodic. 
Thus, any free strongly ergodic aperiodic pmp $\Gamma$-action is in fact OE-superrigid, and 
the proof of Proposition 5.1 applies.  

Let now $\Gamma = \Gamma_1 *_H \Gamma_2$ be as in  $2^\circ$. Recall from (Theorem 1.3 in [K2]) that 
for such a group $\Gamma$, any free ergodic pmp $\Gamma$-action which is aperiodic when restricted to $H$ is OE-superrigid. 
But if $\Gamma \curvearrowright X$ is strongly ergodic and aperiodic on $H$, the above remark shows 
that it is strongly ergodic on any 
finite index subgroup of $H$. Thus, the ultrapower action $\Gamma \curvearrowright A^\omega$ is ergodic on any 
finite index subgroup of $H$, where $A=L^\infty(X)$ as usual.   
Moreover, the same is true for the restriction of this action to any $\Gamma$-invariant 
separable von Neumann subalgebra $A_0\subset A^\omega$ on which $\Gamma$ acts freely. Thus, any such $\Gamma \curvearrowright A_0$ is 
OE-superrigid. Hence, if $\Lambda \curvearrowright Y$ is a free ergodic pmp action that's $\omega$-OE (respectively app-OE) 
with $\Gamma \curvearrowright X$, then the same argument as in the proof of 5.1 applies to get that the two actions follow $\omega$-conjugate 
(resp. app-conjugate). 

Finally, recall from [K1] (respectively [CK]) that if a pmp action $\Gamma \curvearrowright^\sigma X$  is as in part $(a)$ of $3^\circ$ (respectively as 
in part $(b)$ of $3^\circ$) then $\sigma$ is OE-superrigid. But since $\Gamma \curvearrowright A=L^\infty(X)$ is strongly ergodic, by 2.10 so is its 
ultrapower $\Gamma \curvearrowright^{\sigma^\omega} A^\omega$. If in addition $\sigma$ is aperiodic then its 
restriction  to any subgroup $\Gamma_0 \subset \Gamma$ of finite index is still strongly ergodic, so $\Gamma_0 \curvearrowright A^\omega$ is strongly ergodic as well. 
Moreover, the same is true for the restriction of $\sigma^\omega$ to any $\Gamma$-invariant separable subalgebra $A_0\subset A^\omega$ on which $\Gamma$ acts 
freely. Thus, any such $\Gamma \curvearrowright A_0$ follows OE-superrigid. The proof of 5.1 applies again to get the conclusion.   
\hfill $\square$

\proclaim{5.3.  Corollary} Let  $\Gamma=\Gamma_1 \times \Gamma_2$, with $\Gamma_i$ torsion free groups in the class $\Cal C_{reg}$ 
of Monod-Shalom $[\text{\rm MS}]$, $i=1,2$.  Let $\Gamma \curvearrowright X$ be a free pmp $\Gamma$-action whose restriction to $\Gamma_1, \Gamma_2$, 
is strongly ergodic $($e.g., a quotient of a Bernoulli $\Gamma$-action$)$. Then we have: 

$1^\circ$ Any app-OE $($respectively $\omega$-OE$)$ between $\Gamma \curvearrowright X$ and an arbitrary free pmp $\Gamma$-action 
$\Gamma \curvearrowright Y$ comes from an app-conjugacy $($resp. $\omega$-conjugacy$)$. 

$2^\circ$ If in addition $\Gamma_1, \Gamma_2$ have property $(\text{\rm T})$  $($e.g., if each $\Gamma_i$ is an arithmetic lattice in some $Sp(n,1)$, $n\geq 2$, or a quotient  
of such a group$)$, then any app-OE between $\Gamma \curvearrowright X$ and an arbitrary free pmp $\Gamma$-action 
$\Gamma \curvearrowright Y$ comes from a conjugacy. 
\endproclaim
\noindent
{\it Proof}. $1^\circ$ Recall from (Theorem 1.6 in [MS]) that if $\Gamma=\Gamma_1 \times \Gamma_2$ is as in the hypothesis, then an OE between a free pmp 
action $\Gamma \curvearrowright^\sigma X$ whose restriction to $\Gamma_1, \Gamma_2$ is still ergodic, and an arbitrary free pmp action $\Gamma \curvearrowright Y$, 
comes from a conjugacy. 

But if in addition $\sigma_{|\Gamma_i}$ are strongly ergodic, $i=1,2$, then by Proposition 2.10 any restriction 
of the ultrapower action $\Gamma \curvearrowright^{\sigma^\omega} A^\omega$ to a $\Gamma$-invariant  subalgebra $\tilde{A}\subset A^\omega$ 
is still strongly ergodic on $\Gamma_1, \Gamma_2$ (so in particular ergodic). 

Then the same argument as in the proof of 5.1 applies to show that 
any $\omega$-OE (respectively app-OE) between $\Gamma \curvearrowright^\sigma X$ and another free pmp $\Gamma$-action comes from a 
$\omega$-conjugacy (respectively an app-conjugacy).   

Part $2^\circ$ is then an immediate consequence of the first part and 3.11. 
\hfill $\square$

\vskip .1in 
\noindent
{\bf 5.4. Definition.} Let  $A\subset M$ be a Cartan inclusion and $\omega$ a free ultrafilter on $\Bbb N$. We will define the 
$\omega$ and {\it app symmetry groups} of $A\subset M$ as the corresponding (genuine) symmetry groups of the 
Cartan inclusion $A^\omega \subset M(\omega)$, where  $M(\omega)$ denotes, as before, the von Neumann algebra 
generated inside $M^\omega$ by $A^\omega$ and $M$. Thus, if $\Cal U \subset \Cal N_M(A)$ is any group 
with $\Cal U \vee A =M$, then $M(\omega)= \Cal U \vee A^\omega=\overline{\text{\rm sp}\Cal UA^\omega}$. 

It is useful to note that if $B\subset N$ denotes the Cartan inclusion 
$(A\subset M)^t$, then $B^\omega \subset N(\omega)$ is naturally isomorphic to $(A^\omega \subset N(\omega))^t$, 
with the corresponding identification between $(A^t)^\omega$ and $(A^\omega)^t$ being of ultraproduct type.  

\vskip .05in

$1^\circ$ We denote by  Aut$_\omega(A\subset M)$ 
the automorphism group of $A^\omega \subset M(\omega)$, as defined in 2.15.1$^\circ$, and by $\text{\rm Aut}_{app, \omega}(A\subset M)$ 
its  subgroup of automorphisms which restricted to $A^\omega$ are of ultrapower form. We denote $\text{\rm Out}_\omega(A \subset M)$ 
(respectively Out$_{app,\omega}(A\subset M)$) the quotient of Aut$_\omega(A\subset M)$ (resp. of Aut$_{app,\omega}(A\subset M)$) 
by the subgroup of automorphisms implemented by 
unitaries in $\Cal N_{M(\omega)}(A^\omega)$. 

If $\Gamma \curvearrowright^\sigma X$ is a pmp action of a countable group, then we denote Aut$_\omega(\Cal R_\Gamma)$, Out$_\omega(\Cal R_\sigma)$, 
and Aut$_{app,\omega}(\Cal R_\sigma)$, Out$_{app,\omega}(\Cal R_\sigma)$ the corresponding groups of the Cartan inclusion $L^\infty(X)\subset L(\Cal R_\sigma)$. 
Note that Out$_\omega(\Cal R_\sigma)$ (respectively Out$_{app,\omega}(\Cal R_\sigma)$) is an $\omega$-OE (resp. app-OE) invariant for $\sigma$. 

\vskip .05in
$2^\circ$ We denote $Z^1_\omega(A\subset M)=Z^1(A^\omega \subset M(\omega))$, 
$H^1_\omega(A\subset M)=H^1(A^\omega \subset M(\omega))$, 
the first cohomology groups of $A^\omega \subset M(\omega)$, as defined in 2.15.2$^\circ$. We identify these groups, as in 2.15.2$^\circ$,  
with  Aut$_0(A^\omega \subset M(\omega))$,  respectively Out$_0(A^\omega \subset M(\omega))$,  of (classes  
of) automorphisms of $A^\omega \subset M(\omega)$ that leave $A^\omega$ fixed.  
Note that such automorphisms are of ultrapower form when restricted to 
$A^\omega$, so they can be viewed as subgroups of Aut$_{app,\omega}(A\subset M)\subset \text{\rm Aut}_\omega(A\subset M)$. 

If $L^\infty(X)=A \subset M=L^\infty(X)\rtimes_\sigma \Gamma$ for some free ergodic pmp action $\Gamma \curvearrowright^\sigma X$ then,  
with the notations in 2.15.2$^\circ$, we have $Z^1_\omega (A\subset M) = Z^1(\sigma^\omega)$, $H^1_\omega (A\subset M) = H^1(\sigma^\omega)$, 
the first cohomology groups of 
$\Gamma \curvearrowright^{\sigma^\omega} A^\omega$. Thus, an element in 
$Z^1_\omega(A\subset M)$ is given   
by  a map $c:\Gamma \rightarrow \Cal U(A^\omega)$ satisfying $c_g\sigma^\omega_g(c_h)=c_{gh}$, $\forall g,h\in \Gamma$. 
As pointed out in 2.15.2$^\circ$, the topology inherited 
from $\text{\rm Aut}_{\omega}(A\subset A \rtimes_\sigma \Gamma) = 
\text{\rm Aut}(A^\omega \subset A^\omega \rtimes_\sigma \Gamma)$ (of point-$\| \ \|_2$ convergence) corresponds in $Z^1_\omega(\sigma)$ to the point-$\| \ \|_2$-convergence of cocycles, viewed  as functions $\Gamma \rightarrow \Cal U(A^\omega)$.

By writing 
$c_g=(c_{g,n})_n\in \Cal U(A^\omega)$ with $c_{g,n}\in \Cal U(A), \forall n$, we see that $c$ corresponds to a sequence of functions 
$c_n:\Gamma \rightarrow \Cal U(A)$ 
which satisfy the approximate cocycle relation $\lim_\omega \|c_{g,n}\sigma_g(c_{h,n})-c_{gh,n}\|_2=0$, $\forall g,h\in \Gamma$. Two 
such sequences $c=(c_n)_n$, $c'=(c'_n)_n$ give the same element in $Z^1_\omega(\sigma)$ if $\lim_\omega \|c'_{g,n}-c_{g,n}\|_2=0$, $\forall g\in \Gamma$. This shows that 
any sequence of  cocycles $c_n\in Z^1(\sigma)$ gives rise to a $\sigma^\omega$-cocycle $c=(c_n)_n$.  
We denote by $Z^1(\sigma)^\omega$ the subgroup of such cocycles. It is easy to see that it is closed in $Z^1(\sigma^\omega)$. Also, we always have 
$(B^1(\sigma))^\omega= B^1(\sigma^\omega)$. 

By  [Sc1] we have  
$B^1(\sigma)$ closed in $Z^1(\sigma)$  iff $\sigma$ is strongly ergodic and this condition is easily seen to imply that $B^1(\sigma^\omega)$ is closed in $Z^1(\sigma^\omega)$. 
However, the converse is not necessarily true. For instance,  the co-boundary group of $\sigma^\omega$, $B^1(\sigma^\omega)$ (which we saw is equal to $(B^1(\sigma))^\omega$) 
coincides with the closed subgroup  $Z^1(\sigma)^\omega$ of $Z^1(\sigma^\omega)$ whenever $B^1(\sigma)$ is dense in $Z^1(\sigma)$ (e.g., when $\Gamma$ is amenable).  

Note that if $H^1(\sigma)$ is discrete (e.g., if $\Gamma$ has property (T), see [Sc1]), 
then the embedding $(Z^1(\sigma))^\omega \subset Z^1(\sigma^\omega)$ 
induces an embedding of the ultrapower group $H^1(\sigma)^\omega$ into $H^1_\omega(\sigma)$. The interesting problem here is to 
show that, in certain cases, these groups coincide,  $Z^1_\omega(\sigma)=Z^1(\sigma)^\omega$, $H^1_\omega(\sigma)=H^1(\sigma)^\omega$, i.e., 
any approximate cocycle comes from a cocycle. This would show, for instance, that if $H^1(\sigma)$ is finite, 
then $H^1_\omega(\sigma)=H^1(\sigma)$.

\vskip .05in 

3$^\circ$ We denote by $\Cal F_\omega(A\subset M)$ the fundamental group of $A^\omega \subset M(\omega)$ and by $\Cal F_{app,\omega}(A\subset M)$  
the subgroup of  all $t>0$ with the property that there exists an isomorphism $(A^\omega \subset M(\omega)) \simeq (A^\omega \subset M(\omega))^t$ which 
is of ultraproduct form when restricted to $A^\omega$. Note that we have natural embeddings $\Cal F(A\subset M)\subset \Cal F_{app}(A\subset M) 
\subset \Cal F_\omega(A\subset M)$ and that the group $\Cal F_\omega(A\subset M)$ is an $\omega$-OE invariant while $\Cal F_{app}(A\subset M)$ is an 
app-OE invariant.  

\vskip .1in

By 3.11.2$^\circ$, 3.12.2$^\circ$, Corollary 5.3 above and (7.13.3 in [IPP]), we thus obtain the following concrete calculations of  Out$_{app,\omega}(A\subset M)$ 
for group measure space Cartan inclusions coming from some special classes of actions of property (T) groups. 

\proclaim{5.5. Corollary}  Let $\Gamma$ be an ICC hyperbolic  group with the property $(\text{\rm T})$ and $\Gamma \times \Gamma \curvearrowright^{\sigma_i} X_i^\Gamma$ 
be double Bernoulli $\Gamma \times \Gamma$-actions with finite base $(X_i, \mu_i)$, $i\in I$. 
Two such actions $\sigma_i, \sigma_j$ are app-OE if and only if $(X_i, \mu_i)\simeq (X_j, \mu_j)$. 
Also, $\text{\rm Out}_{app,\omega}(\Cal R_{\sigma_i})=\text{\rm Aut}(X_i, \mu_i)$.  
\endproclaim

On the other hand, by Proposition 2.12.3$^\circ$ and the remark at the end of 2.15.3$^\circ$, we  get:  

\proclaim{5.6. Corollary}  If $\Gamma$ has one non-zero, finite $L^2$-Betti number, then 
any free ergodic pmp $\Gamma$-action  has trivial $\omega$-fundamental group, $\Cal F_\omega(\Cal R_\Gamma)=\{1\}$.  

\endproclaim

\vskip .1in 
\noindent
{\bf 5.7. Final remarks}. It is tempting to re-examine existing OE rigidity results and calculations of invariants 
and try to use them, or try to adapt their proofs, to obtain app-OE rigidity results. But the only OE-rigidity  results that 
we were able to use to derive app-OE/$\omega$-OE rigidity were the ones in [K1, K2], [CK], [MS]. This is because the OE-rigidity results 
in these papers have almost no restrictions on the ``source action'' (albeit putting much restrictions on the acting group). 

In turn, the OE rigidity results that require very specific restrictions on the source action, like being Bernoulli, Gaussian, or profinite, cannot be applied 
directly, and the deformation-rigidity techniques used to prove them do not seem  to adapt well to the ``ultrapower framework''
(e.g., because separability and deformation arguments fail). 

So a general problem that's of interest in this respect is to search for more examples of OE-superrigid group actions $\Gamma \curvearrowright^\sigma X$ 
where the only conditions that  would be imposed on the action $\sigma$ have good  ``$\omega$-permanence'', i.e., they ``survive'' 
to all separable sub-actions of the ultrapower $\Gamma$-action $\sigma^\omega$. For once one has such a result, the proof of 
Proposition 5.1 applies and we can derive app-OE rigidity for $\sigma$. 

Such OE-rigidity results are particularly interesting to have for property (T) groups $\Gamma$, because as we have seen in Corollary 3.11, 
app-conjugacy for such actions is the same as conjugacy. This can also lead to complete calculations of app-symmetry groups, like in Corollary 5.5 above.

\vskip .05in 
\noindent
{\it 5.7.1. Calculations of $\text{\rm Out}_{app,\omega}(\Cal R_\Gamma)$}. Strong rigidity results like Corollaries 5.3, 
or like Corollary 5.2 applied to just the case $\Gamma = \Lambda$, reduce the calculation of 
the app-Outer automorphism group $\text{\rm Out}_{app,\omega}(\Cal R_\Gamma)$ to the calculation of the app-conjugacy invariant 
Aut$^\sigma_{app,\omega}(X)$, defined in 3.11.2$^\circ$ as the group of sequences of automorphisms $\theta_n$ of $(X, \mu)$ 
that $\omega$-asymptotically commute with $\sigma_g$, $\forall g $. To see this, we send the reader to page 445 in [P7], or to similar considerations in [F2], [MS]   
(we have already seen in 5.5 above a sample of how this works). 

So in order to produce more classes of group actions $\Gamma \curvearrowright^\sigma X$ with calculable $\text{\rm Out}_{app,\omega}(\Cal R_\sigma)$ 
one needs to obtain more classes of groups $\Gamma$ and actions $\Gamma \curvearrowright X$ with the property that,   
for any $\Gamma$-invariant separable subalgebra $A_0\subset A^\omega$ that contains $A$, any automorphism of $A_0\subset A_0 \rtimes \Gamma$ 
comes from a self-conjugacy of $\Gamma \curvearrowright A_0$.

\vskip .05in 
\noindent
{\it 5.7.2. Untwisting app-cocycles.}  We leave wide open the problem of concrete calculations 
of the app-OE invariant $H^1_\omega(\sigma)$ of a free ergodic pmp action $\Gamma \curvearrowright ^\sigma X$.  As explained before, 
the actions to consider should be strongly ergodic, a property that's automatic when $\Gamma$ has property (T).  
For such group actions, the result to seek is that the natural embedding $(H^1(\sigma))^\omega \subset H^1_\omega(\sigma)$ 
is in fact an equality. 

This amounts to proving that if $c_n: \Gamma \rightarrow \Cal U(A)$ are maps satisfying 
$$
\lim_{n \rightarrow \omega} \|c_n(g)\sigma_g(c_n(h))-c_n(gh)\|_2=0, \forall g, h\in \Gamma,
$$ 
then there exist $c'_n \in Z^1(\sigma)$ such that $\lim_\omega \|c'_n(g)-c_n(g)\|_2 = 0$, $\forall g$. Once this is proved, 
if $H^1(\sigma)$ itself is calculable (e.g., when $\sigma$ belongs to calculable cohomology in [Ge], [PS], [P6]), we are done. 
Or we can try to prove directly that the cocycle $c=(c_n)_n$ for $\Gamma \curvearrowright^{\sigma^\omega} A^\omega$ ``untwists''   (like in [Ge],  [PS], [P6]). 

We can in fact view this problem as a particular case of a general question 
about untwisting app-cocycle  with $\Cal U_{fin}$ targets, in the spirit of ([P8]). While this is somewhat long and tedious to formulate 
in its full generality, we will only state it in the case of discrete targets (the other interesting case being scalar targets, which we have already stated above).

Thus, let $\Gamma \curvearrowright^\sigma (X,\mu)$ be a free ergodic pmp action of a countable group and  
$\Lambda$ another countable group. We let $\Gamma$ act on the space $\Lambda^X$ of 
measurable maps from $X$ to $\Lambda$ (viewed as a 
subgroup of unitaries in $L^\infty(X)\overline{\otimes} L(\Lambda)$) 
by $\sigma_g(b)(t)=b(\sigma_{g^{-1}}(t))$, $\forall t\in X$, $g\in \Gamma$. 
An app-cocycle (respectively $\omega$-cocycle) for $\sigma$ 
with target group $\Lambda$ is a sequence of maps $c_n: \Gamma \rightarrow \Lambda^X$ with the property that 
$\lim_{n \rightarrow \infty} \|c_n(g)\sigma_g(c_n(h)) - c_n(gh)\|_2=0$ (resp. $\lim_{n \rightarrow \omega} \|c_n(g)\sigma_g(c_n(h)) - c_n(gh)\|_2=0$), $\forall g, h \in \Gamma,$
and such that for each $g\in \Gamma$, $h\in \Lambda$, we have $\mu(\{t\in X \mid c_n(g)(t)=h\})$ constant in  $n$. 

We say that $c=(c_n)_n$ can be app-untwisted (resp. $\omega$-untwisted) if there exist a group morphism $\delta: \Gamma \rightarrow \Lambda$ and a 
sequence of maps $w_n \in \Lambda^X$ 
such that $\lim_{n\rightarrow \infty} \| w_n c_n(g)\sigma_g(w_n^*) - \delta(g)\|_2=0$  
(resp. $\lim_{n\rightarrow \omega} \| w_n c_n(g)\sigma_g(w_n^*) - \delta(g)\|_2=0$), $\forall g$, and such that 
for each $h\in \Lambda$, $\mu(\{t\in X \mid w_n(t)=h\})$ is constant in $n$.

We say that $\Gamma \curvearrowright^\sigma X$ is {\it app-cocycle superrigid} (resp. $\omega$-{\it cocycle superrigid}) 
if any app-cocycle (resp. $\omega$-cocycle) with discrete targets can be app-untwisted 
(resp. $\omega$-untwisted).  Arguing like on (page 287 of [P8]), 
we see that any approximate orbit equivalence $\theta=(\theta_n)_n$ 
(resp. $\omega$-OE $\theta$) between two 
free ergodic pmp actions $\Gamma \curvearrowright^\sigma X$, $\Lambda \curvearrowright^\rho Y$  gives rise to an 
app-cocycle (resp. $\omega$ cocycle) $c=(c_n)$, with $c_n:\Gamma \rightarrow Y^\Lambda$, by writing each $\theta \sigma^\omega(g) \theta^{-1}$ 
as $(\oplus_h \rho(h)_{|q^g_{n,h}B})_n$ and letting $c_n(g)$ be equal to $h$ when restricted to the subset of $Y$ 
corresponding to  $q^g_{n,h}$. As in ([P8]), untwisting the app-cocycle $c=(c_n)_n$ by some $w=(w_n)_n$ in the above sense   
means that $w_n\in \Lambda^Y$ 
can be interpreted as elements in the normalizer of $B^\omega$ in $N(\omega)$ satisfying 
$w\theta \sigma^\omega(\Gamma)\theta^{-1}w^*
\subset \rho^\omega(\Lambda)$, with the fact that $\tilde{\theta} (M(\omega))=N(\omega)$ implying that we actually have equality (where $\tilde{\theta}$ is the extension of 
$\theta:L^\infty(X)^\omega\simeq L^\infty(Y)^\omega$ given by 4.1).  
In other words, $\theta$ and $w$ implement an app-conjugacy (resp. $\omega$-conjugacy) of the two actions.  

Thus, if we can prove app-cocycle superrigidity 
(resp. $\omega$-cocycle superrigidity) for a free pmp $\Gamma$ action $\sigma$, then we also get app-OE superrigidity 
(resp. $\omega$-OE superrigidity)  for $\sigma$. 

\vskip .05in 
\noindent
{\it 5.7.3. On the existence of non weak-OE actions}. As we mentioned before, we leave open the 
question of whether  every property (T) group $\Gamma$ has two non weak-OE free ergodic pmp actions. 
And further, whether any non-amenable group $\Gamma$ has ``many'' non weak-OE group actions, or at least 
that it has ``many'' non app-OE actions (see [AE] for partial answers to this problem). 

For the ``classic'' OE-case of this question, the successful final answer in ([I1], [GL], [E]) used 
deformation-rigidity arguments, by exploiting the existence of many 
free ergodic pmp $\Bbb F_2$-actions with relative property (T) in [GP]. 
But it seems difficult to  make such arguments work in ultrapower framework. In the absence of that, the calculation 
of app-cohomology $H^1_\omega(\sigma)$ as described in 5.7.2 would at least solve the problem for 
existence of many non app-OE actions of property (T) groups. 

For the weak-OE version of this problem, another approach is to seek properties of an action that have good $\omega$-permanence. 
In this respect, it would be interesting to investigate whether the Cartan inclusion $A^\omega \subset M(\omega)=A^\omega \rtimes_{\sigma^\omega} \Gamma$ 
associated with  a Bernoulli $\Gamma$-action $\sigma$ of a non-amenable group  
$\Gamma$ can admit a sub- Cartan inclusion $B\subset N$ that has relative property (T) (or is rigid) in the sense of ([P5]), or at least that there exists 
$A_0\subset A^\omega$ diffuse such that $A_0\subset M(\omega)$ is rigid. Another question is whether $A^\omega \subset M(\omega)$ 
can admit sub-Cartan inclusions $B\subset N$ that come from profinite actions (e.g., $B\subset N$ weakly-compact, in the sense  of [OP]).

\vskip .05in 
\noindent
{\it 5.7.4.  A connection to $2$nd cohomology}. A well known open problem in OE rigidity is to effectively calculate the 2-cohomology group 
$H^2(\sigma)$ of a free ergodic pmp action $\Gamma \curvearrowright^\sigma X$ (see 6.6 in [P8]). Recall from Example 2.4  that this (abelian) group is defined as the quotient between  
the group of (normalized) 2-cocycles for $\sigma$, 
$$
Z^2(\sigma)=\{ v: \Gamma \times \Gamma \rightarrow \Cal U(A) \mid 
v_{g,h}v_{gh,k}=\sigma_g(v_{h,k})v_{g,hk},   \forall g,h,k \}, 
$$ 
by the subgroup of (normalized) co-boundary cocycles 
$$
B^2(\sigma)=\{ v\in Z^2(\sigma) \mid \exists w:\Gamma \rightarrow \Cal U(A), v_{g,h}=w_g \sigma_g(w_h)w_{gh}^*, \forall g,h\}, 
$$

We endow the group $Z^2(\sigma)$  with the topology of  pointwise convergence in the norm $\|  \ \|_2$. It is a Polish 
group with respect to this topology.  

Note that $B^2(\sigma)$ is the image under the {\it boundary map}  $\partial$ from the product group $\Cal U(A)^{\Gamma}$, 
defined by $\partial w (g,h)=w_g \sigma_g(w_h) w_{gh}^*$,  $\forall g,h\in \Gamma$. 
The group morphism $\partial$ is clearly continuous and its kernel is the (closed) 
subgroup $Z^1(\sigma)$ of $1$-cocycles for $\sigma$. Thus, $\partial$ induces a group isomorphism $\tilde{\partial}$ 
from the (Polish) group $\Cal U(A)^{\Gamma}/Z^1(\sigma)$ onto the (topological) group $B^2(\sigma)$. But while $\tilde{\partial}$ 
is continuous, its inverse is continuous iff $B^2(\sigma)$ is closed in $Z^2(\sigma)$, see Proposition 5.8 below.

It is a consequence 
of [CFW] that $H^2(\sigma)=\{1\}$ for any ergodic pmp action $\sigma$ of an amenable group $\Gamma$. 
But surprisingly, as of now, there has been no calculation of $H^2(\sigma)$ for a free ergodic pmp action of a non-amenable group! 
As already mentioned in Example 2.4, any element in the 2nd cohomology group of $\Gamma$ 
$$
Z^2(\Gamma, \Bbb T)=\{\lambda: \Gamma \times \Gamma \rightarrow \Bbb T \mid 
\lambda_{g,h}\lambda_{gh,k}=\lambda_{h,k}\lambda_{g,hk}, \forall g,h,k\in \Gamma \}
$$
implements an element in $Z^2(\sigma)$. But it is not easy to decide when a scalar cocycle gives a non-trivial  
element in $H^2(\sigma)$.  One situation that's particularly desirable is for $H^2(\sigma)$ to actually coincide with $H^2(\Gamma, \Bbb T)$, 
something we believe should be true  in certain cases (e.g., when $\Gamma$ has some strong form of property (T) and $\sigma$ is Bernoulli). 
For  this to happen (and more generally for $H^2(\sigma)$ to be calculable) a necessary condition is that 
$H^2(\sigma)$ be a Polish group, i.e., $B^2(\sigma)$  be closed 
in $Z^2(\sigma)$. 

The next Proposition shows that $B^2(\sigma)$ is closed iff the  condition $(Z^1(\sigma))^\omega=Z^1(\sigma^\omega)$, that we 
discussed in 5.7.2 and which is needed in the calculation of the app 1-cohomology group $H^1_\omega(\sigma)$, holds true. This partly explains why  
both the calculation of $H^1_\omega(\sigma)$ and $H^2(\sigma)$ are so difficult.

\proclaim{5.8.  Proposition} The following conditions 
are equivalent: 
\vskip .05in
$(a)$ $Z^1(\sigma)^\omega = Z^1(\sigma^\omega)$; 
\vskip .05in 
$(b)$ $B^2(\sigma)$ is closed in $Z^2(\sigma)$; 
\vskip .05in 
$(c)$ The boundary map $\partial: \Cal U(A)^{\Gamma} \rightarrow B^2(\sigma)$ 
is open. In other words, the group isomorphism $\Cal U(A)^{\Gamma}/Z^1(\sigma)\simeq B^2(\sigma)$,  
implemented by the boundary map, is an isomorphism of Polish groups.

\endproclaim 
\noindent
{\it Proof}. The proof follows closely arguments in [C2]. We will prove that (a) $\Rightarrow$ (c) $\Rightarrow$ (b) $\Rightarrow$ (a). 

(a) $\Rightarrow$ (c): To show that $\partial$ is open, it is sufficient to prove  
that if a sequence of 2-coboundaries $\tilde{v}_{g,h}^n = \tilde{w}_g^n\sigma_g(\tilde{w}_h^n)(\tilde{w}_{gh}^n)^*$ converges to $1$, then there is a
subsequence $(\tilde{w}^{n_k})$ and $s^k\in Z^1(\sigma)$ such that $s^k_g \tilde{w}^{n_k}_g\to 1$ strongly, $\forall g\in \Gamma$.
So suppose that $\tilde{w}_g^n\sigma_g(\tilde{w}_h^n)(\tilde{w}_{gh}^n)^*\to 1$ strongly, $\forall g,h\in \Gamma$. Then $g\mapsto (\tilde{w}^n_g)_n$ belongs to
$Z^1(\sigma^\omega)$. Since $Z^1(\sigma)^\omega = Z^1(\sigma^\omega)$, there is a sequence $(w_g^n)_{g\in G}$ of 1-cocycles for $\sigma$ such that
$\lim_{n\to\omega}\|w_g^n-\tilde{w}_g^n\|_2 = 0$, $\forall g\in \Gamma$. We claim that this implies the desired conclusion.
Indeed, as $\Gamma$ is countable, we can construct an increasing map $\Bbb N\ni k\mapsto n_k\in\Bbb N$ such that $\|\tilde{w}_g^{n_k}-w_g^{n_k}\|_2\to 0$ 
as $k\to\infty$, $\forall g\in \Gamma$, and then $((w^{n_k})^{-1}\tilde{w}^{n_k})$ converges to $1$ in the topology on $\Cal U(A)^{\Gamma}$.

(c) $\Rightarrow$ (b): View $\Cal U(A)^{\Gamma}$ as a subset of 
			$N= \bigoplus_{g\in \Gamma} (A)_g = 
			\ell^\infty(\Gamma)\bar{\otimes} A$. Put $V_n = \{a\in 
			(N)_1\,:\,\|a\xi_j\|<2^{-n}\text{ for 
			}1\leq j\leq n\}$, where $\{\xi_j\}_{j=1}^\infty$ is an orthonormal basis in 
			$\ell^2(\Gamma)\bar{\otimes}L^2(A)$. Note that $w^n\to 0$ strongly in $N = 
			\ell^\infty(\Gamma)\bar{\otimes} A$ if and only if $(w^n_g)_{g\in \Gamma}\to 
			0$ in 
			$\Cal U(A)^\Gamma$ and that the $V_n$ form a decreasing sequence of open neighborhoods of
			$0\in N$ (in the 
			strong operator topology restricted to the unit ball) that shrinks 
			to a point and has 
			the property that $uV_n\subset V_n$ for all $u\in\Cal U(A)^\Gamma$. Choose 
			open 
			sets $W_n\subseteq 
			\partial((V_n+Z^1(\sigma))\cap \Cal U(A)^\Gamma)$ that contain $1$, 
			noting 
			that 
			$(V_n+Z^1(\sigma))\cap \Cal U(A)^\Gamma$ is an open set. Then
			$$
				\partial(u)\in W_n\Rightarrow u\in V_n+Z^1(\sigma).
			$$
			Let $v\in \overline{B^2(\sigma)}$ be given. Choose 
			$\partial(w^n)\in 
			B^2(\sigma)$ such that $\partial(w^n)\to v$ and 
			$\partial((w^{n+1})^*w^n)\in W_n$ for all $n$. Then 
			$(w^{n+1})^*w^n\in 
			V_n+Z^1(\sigma)$, i.e., 
			$w^{n} = a_n + w^{n+1}s^n$ with $a_n\in V_n$, $s^n\in Z^1(\sigma)$. So by 
			multiplying 
			each $w^n$ by a 
			1-cocycle if necessary, we may arrange that $(w^n_g)_{g\in \Gamma}$ has 
			a strong limit 
			in $N$ (as $(w^n\xi_j)_n$ is Cauchy for each $j$, converging to 
			$\eta_j$, say, and we can define an operator $w$ by 
			$w(\sum\alpha_j\xi_j) 
			= 
			\sum\alpha_j\eta_j$, which is certainly the pointwise limit of the 
			$(w^n)$
			on the dense linear span of the $\xi_j$, hence is the strong limit 
			of 
			the $(w^n)$), 
			and therefore also that $(w^n_g)$ has 
			a
			strong limit $w_g\in \Cal U(A)$, $\forall g\in \Gamma$ (defined by $w_g(\xi) = 
			w(\delta_g\otimes \xi)$, $\forall g\in \Gamma$), which is necessarily 
			unitary, since $A$ is abelian. 
			It follows that $v = \partial(w)\in B^2(\sigma)$.
			
(b) $\Rightarrow$ (a): Suppose that $B^2(\sigma)$ is 
			a closed subgroup of 
					$Z^2(\sigma)$. Then $B^2(\sigma)$ is a Polish space,
					so the continuous bijective homomorphism
					$$
						\tilde{\partial}\colon\Cal U(A)^\Gamma/Z^1(\sigma)\to 
						B^2(\sigma)
					$$
					automatically has a continuous inverse map (cf.\ e.g.\ a 
					lemma in 
					[C2]). Let $g\mapsto 
					(\tilde{w}_g^n)$ be an element of $Z^1(\sigma^\omega)$. Then
					$$
							\lim_{n\to\omega} 
									\|\partial(\tilde{w}^n)_{g,h}-1\|_2=\lim_{n\to\omega}
						\|\tilde{w}_g^n\sigma_g(\tilde{w}_h)(\tilde{w}_{gh}^n)^*-1\|_2
						 = 0
					$$
					for all $g,h\in \Gamma$. Since $\Cal U(A)^{\Gamma\times \Gamma}$ is equipped 
					with the product topology, we infer that 
					$\lim_{n\to\omega}\tilde{\partial}(\overline{\tilde{w}^n}) 
					= 1$, where $\overline{\tilde{w}^n}$ is the class of $\tilde{w}^n$ in $\Cal U(A)^\Gamma/Z^1(\sigma)$. By continuity of $\tilde{\partial}^{-1}$, it follows easily
					that $\lim_{n\to\omega}  \overline{\tilde{w}^n} = 
					\overline{1}$, hence that $(\tilde{w}_g^n)\in 
							 Z^1(\sigma)^\omega$.

\hfill $\square$

\head  References \endhead

\item{[AE]} M. Abert, G. Elek: {\it Dynamical properties of profinite actions},  Ergodic Theory Dynam. Systems {\bf 32} (2012), 1805-1835.

\item{[A]} C. Anantharaman-Delaroche: {\it On Connes' property} (T) {\it for von Neumann algebras}, Math. Japon. {\bf 32} (1987), 337-355.

\item{[AW]} M. Abert, B. Weiss: {\it Bernoulli actions are weakly contained in any free action},  Ergodic Theory Dynam. Systems {\bf 33} (2013), 323-333.

\item{[Bo]} F. Boca: {\it On a method for constructing irreducible finite index subfactors of Popa}, Pacific J. Math. {\bf 161} (1993), 201-231. 

\item{[B1]} L. Bowen: {\it Measure conjugacy invariants for actions of countable sofic groups}, J. Amer. Math. Soc. {\bf 23} (2010), 217-245.

\item{[B2]} L. Bowen: {\it Weak density of orbit equivalence classes of free group actions}, To appear in Groups, Geometry and Dynamics.

\item{[BHI]} L. Bowen, D. Hoff, A. Ioana: {\it von Neumann's problem and extensions of non-amenable equivalence relations}, 
arXiv:1509.01723

\item{[CK]} I. Chifan, Y. Kida: {\it OE and W$^*$ superrigidity results for actions by surface braid groups}, arXiv:1502.02391

\item{[CKT]} C. Conley, A. Kechris and R. Tucker-Drob: {\it Ultraproducts of measure preserving actions and graph combinatorics}, 
Erg. Th. and Dyn. Syst., {\bf 33} (2013), 334-374.

\item{[C1]} A. Connes: {\it Classification of injective factors}, Ann. of Math., {\bf 104} (1976), 73-115.

\item{[C2]} A. Connes: {\it Almost periodic states and factors of type III$_1$}, J. Funct. Anal., {\bf 16} (1974), 415-445.

\item{[CFW]} A. Connes, J. Feldman, B. Weiss: {\it An amenable equivalence relation is generated by a single
transformation}, Erg. Theory Dyn. Sys.  {\bf 1} (1981), 431-450.

\item{[CJ]} A. Connes, V.F.R. Jones: {\it A} II$_1$ {\it factor with two non-conjugate Cartan subalgebras}, 
Bull. Amer. Math. Soc. {\bf 6} (1982), 211Ð212.

\item{[CW]} A. Connes, B. Weiss: {\it Property} (T) {\it and
asymptotically invariant sequences}, Israel J. Math. {\bf 37} (1980), 209-210. 

\item{[D]} J. Dixmier: {\it Sous-anneaux ab\' eliens maximaux dans les facteurs de type fini}, Ann. of Math. {\bf 59} (1954), 279-286. 

\item{[Dy]}  H. Dye: {\it On groups of measure preserving
transformations} I, Amer. J. Math, {\bf 81} (1959), 119-159. 

\item{[EL]} G. Elek, G. Lippner: {\it Sofic equivalence relations}, JFA {\bf 258} (2010) 1692-1708.

\item{[E]} I. Epstein: {\it Orbit inequivalent actions of non-amenable groups}, arXiv:0707.4215 

\item{[FM]} J. Feldman, C.C. Moore: {\it Ergodic equivalence
relations, cohomology, and von Neumann algebras} I, II, Trans. AMS {\bf 234} (1977), 289-324, 325-359. 

\item{[Fu1]} A. Furman:  {\it Orbit equivalence rigidity}, Ann. of Math. {\bf 150} (1999), 1083-1108.

\item{[Fu2]} A. Furman: {\it Outer automorphism group of some equivalence relations}, Comment. Math. Helv. {\bf 80} (2005), 157-196. 

\item{[G1]} D. Gaboriau: {\it Cout des r\'elations
d'\'equivalence et des groupes}, Invent. Math. {\bf 139} (2000),
41-98.

\item{[G2]} D. Gaboriau: {\it Invariants $\ell^2$ de r\'elations
d'\'equivalence et de groupes},  Publ. Math. I.H.\'E.S. {\bf 95}
(2002), 93-150.

\item{[GL]} D. Gaboriau, R. Lyons: {\it A Measurable-Group-Theoretic Solution to von Neumann's Problem}, 
Invent. Math., {\bf 177} (2009), 533-540.

\item{[GP]}  D. Gaboriau, S. Popa: {\it An uncountable family of non-orbit equivalent actions of $\Bbb F_n$} J. Amer. Math. Soc. {\bf 18} (2005), 547-559.

\item{[Ge87]} S.L. Gefter: {\it On cohomologies of ergodic actions
of a T-group on homogeneous spaces of a compact Lie group}
(Russian), in ``Operators in functional spaces and questions of
function theory'', Collect. Sci. Works, Kiev, 1987, pp 77-83.

\item{[H]} G. Hjorth: {\it A converse to Dye's Theorem}, Trans. AMS. {\bf 357} (2005), 3083-3103.  

\item{[Ho]} C. Houdayer: {\it Invariant percolation and measured theory of nonamenable groups} (after Gaboriau-Lyons, Ioana, Epstein). 
S\' eminaire Bourbaki. Vol. 2010/2011. Ast\' erisque No. {\bf 348} (2012), Exp. No. 1039, 339-374.

\item{[I1]} A. Ioana: {\it Orbit inequivalent actions for groups containing a copy of $\Bbb F_2$}, Invent. Math. {\bf 185} (2011), 55-73.

\item{[I2]} A. Ioana: {\it Cocycle superrigidity for profinite actions of property} (T) {\it groups}, Duke Math J. {\bf 157}, (2011), 337-367.

\item{[I3]} A. Ioana: {\it Classification and rigidity for von Neumann algebras}, European Congress of Mathematics, EMS (2013), 601-625.

\item{[IPP]} A. Ioana, J. Peterson, S. Popa: {\it Amalgamated Free Products of weakly 
rigid Factors and Calculation of their Symmetry Groups},
Acta Math. {\bf 200} (2008), No. 1, 85-153. 

\item{[IT]} A. Ioana, R. Tucker-Drob: {\it Weak containment rigidity for distal actions}, preprint \newline arXiv:1507.05357. 

\item{[Ke]} A. Kechris: ``Global aspects of ergodic group actions'', Mathematical Surveys and Monographs, {\bf 160}. American
Mathematical Society, Providence, RI, 2010.

\item{[K1]} Y. Kida: {\it Orbit equivalence rigidity for ergodic actions of the mapping class group},  Geometria Dedicatae {\bf 131} (2008), 99-109. 

\item{[K2]} T. Kida: {\it Rigidity of amalgamated free products in measure equivalence theory}. J. Topol., {\bf 4} (2011), 687-735. 

\item{[L]} G. Levitt {\it On the cost of generating an equivalence relation}, Ergodic Theory Dynam. Systems {\bf 15} (1995), 1173-1181. 

\item{[Lu]} W. Lueck: {\it The type of the classifying space for a family of subgroups}, J. Pure Appl. Algebra {\bf 149} (2000), 177-203. 

\item{[MS]} N. Monod, Y. Shalom: {\it Orbit equivalence
rigidity and bounded cohomology},  Ann. of
Math. {\bf 164} (2006), 825-878.

\item{[MvN]} F. Murray, J. von Neumann:
{\it On rings of operators}, Ann. Math. {\bf 37} (1936), 116-229.

\item{[OW1]} D. Ornstein, B. Weiss: {\it Ergodic theory of
amenable group actions I. The Rohlin Lemma} Bull. A.M.S. (1) {\bf 2}
(1980), 161-164.

\item{[OW2]} D. Ornstein, B. Weiss: {\it Entropy and isomorphism theorems for actions of amena- \newline ble groups}, Journal d'Analyse Math. {\bf 48} 
(1987), 1-143.

\item{[OP]} N. Ozawa, S. Popa: {\it On a class of} II$_1$ {\it factors with at most one Cartan subalgebra}, Ann. of Math.
{\bf 172} (2010), 101-137.

\item{[PP]} M. Pimsner, S. Popa: {\it Entropy and index for
subfactors}, Annales Scient. Ecole Norm. Sup., {\bf 19} (1986),
57-106.

\item{[P1]} S. Popa: {\it Maximal injective subalgebras in factors associated with free groups}, Adv. Math. {\bf 50} (1983)
27-48.

\item{[P2]} S. Popa: {\it Notes on Cartan subalgebras in type} II$_1$ {\it factors}, 
Mathematica Scandinavica, {\bf 57} (1985), 171-188

\item{[P3]} S. Popa: {\it Correspondences}, INCREST Preprint, 56/1986. 

\item{[P4]} S. Popa: {\it Some properties of the symmetric enveloping algebras with applications to amenability and property} T, 
Documenta Mathematica, {\bf 4} (1999), 665-744.

\item{[P5]}  S. Popa: {\it On a class of type} II$_1$ {\it factors with
Betti numbers invariants}, Ann. of Math {\bf 163} (2006), 809-899
(math.OA/0209310).

\item{[P6]} S. Popa: {\it Some computations of 1-cohomology groups and construction of non-orbit-equivalent actions}, J. Inst. Math. Jussieu {\bf 5} (2006), 309-332.

\item{[P7]} S. Popa: {\it Strong Rigidity of} II$_1$ {\it Factors Arising from Malleable Actions of w-Rigid Groups} I, II, 
Invent. Math., {\bf 165} (2006), 369-408 and 409-453.

\item{[P8]} S. Popa: {\it Cocycle and orbit equivalence superrigidity
for malleable actions of $w$-rigid groups}, Invent. Math. {\bf 170} (2007), 243-295.

\item{[P9]} S. Popa: {\it On the superrigidity of malleable
actions with spectral gap},  J. Amer. Math. Soc. {\bf 21}
(2008), 981-1000 (math.GR/0608429).

\item{[P10]} S. Popa: {\it Deformation and rigidity for group actions and von Neumann algebras}, 
in ÒProceedings of the International Congress of MathematiciansÓ (Madrid 2006), Volume I, EMS Publishing House, Zurich 2006/2007, pp. 445-479.

\item{[P11]} S. Popa: {\it Independence properties in subalgebras of ultraproduct} II$_1$ {\it factors}, 
J. Funct. Anal. {\bf 266} (2014), 5818-5846.

\item{[P12]} S. Popa: {\it A} II$_1$ {\it factor approach to the Kadison-Singer problem}, Comm. Math. Physics. {\bf 332} (2014), 379-414.

\item{[PS]} S. Popa, R. Sasyk:  {\it On the cohomology of Bernoulli actions}, Ergodic Theory Dynam. Systems {\bf 27} (2007), 241-251 
(math.OA/0310211).

\item{[PSV]} S. Popa, D. Shlyakhtenko, S. Vaes: {\it Cohomology and $L^2$-Betti numbers for subfactors and quasi-regular inclusions}, math.OA/1511.07329

\item{[PV]} S. Popa, S. Vaes: {\it On the fundamental group of} II$_1$ {\it factors and equivalence relations arising from group actions}, 
in ÒNon-Commutative GeometryÓ, Proc. Conference in Honor of Alain Connes'  60th birthday, April 2-6, 2007, IHP Paris. 

\item{[Sc1]} K. Schmidt: {\it Asymptotically invariant sequences
and an action of $SL(2, \Bbb Z)$ on the $2$-sphere}, Israel. J.
Math. {\bf 37} (1980), 193-208.

\item{[Sc2]} K. Schmidt: {\it Amenability, Kazhdan's property T, strong ergodicity and invariant means for ergodic group-actions}. 
Erg. Theory Dyn. Systems {\bf 1} (1981), 223-236.

\item{[S]} I.M. Singer: {\it Automorphisms of finite factors},
Amer. J. Math. {\bf 177} (1955), 117-133.

\item{[T]} R. Tucker-Drob: {\it Weak equivalence and non-classifiability of measure preserving actions},  
Ergodic Theory Dynam. Systems {\bf 35} (2015), 293-336.

\item{[V1]} S. Vaes:  {\it Rigidity results for Bernoulli actions and their von Neumann algebras} (after Sorin Popa)
SŽminaire Bourbaki, exposŽ 961. AstŽrisque {\bf 311} (2007), 237-294.

\item{[V2]} S. Vaes: {\it Rigidity for von Neumann algebras and their invariants}, 
In ``Proceedings ICM 2010'', Vol. III, Hindustan Book Agency, 2010, pp. 1624-1650

\item{[W]} F. Wright: {\it Reduction for algebras of finite type}, Ann.  Math. {\bf 60} (1954), 560-570.

\item{[Z]} R. Zimmer: {\it Strong rigidity for ergodic actions of semisimple Lie groups},  Ann. of
Math. {\bf 112} (1980), 511-529.

\enddocument